%% file: steklov_exist.tex
\documentclass[a4paper,draft]{article}
\usepackage{a4wide}
\input{commands.tex}
\usepackage{fancyhdr}
\title{On the existence and stability of modified Maxwell Steklov eigenvalues%
\thanks{The author of this work was supported by DFG project 468728622 and DFG SFB 1456 project 432680300.}
}
\author{Martin Halla%
\thanks{%
Institut f\"ur Numerische und Angewandte Mathematik,
Georg-Augst Universität Göttingen,
Lotzestr.\ 16-18, 37083 Göttingen, Deutschland.
e-mail: m.halla@math.uni-goettingen.de}
}

\begin{document}

\maketitle
\begin{abstract}
\noindent
In recent decades qualitative inverse scattering methods with eigenvalues as target signatures received much attention.
To understand those methods a knowledge on the properties of the related eigenvalue problems is essential.
However, even the existence of eigenvalues for such (nonselfadjoint) problems is a challenging question and existing results for absorbing media are usually established under unrealistic assumptions or a smoothing of the eigenvalue problem.
We present a technique to prove the existence of infinitely many eigenvalues for such problems under realistic assumptions.
In particular we consider the class of scalar and modified Maxwell nonselfadjoint Steklov eigenvalue problems.
In addition, we present stability results for the eigenvalues with respect to changes in the material parameters.
In distinction to existing results the analysis of the present article requires only minimal regularity assumptions.
By that we mean that the regularity of the domain is not required to be better than Lipschitz, and the material coefficients are only assumed to be piece-wise $W^{1,\infty}$.
Also the stability estimates for eigenvalues are obtained solely in $L^p$-norms ($p<\infty$) of the material perturbations.
\\

\noindent
\textbf{MSC:} 35P05, 35R05, 35R30, 78M35
\\

\noindent
\textbf{Keywords:} Steklov eigenvalues, Maxwell's equation, inverse scattering, existence of eigenvalues, stability of eigenvalues
\end{abstract}

\section{Introduction}\label{sec:introduction}

An important unresolved problem in inverse scattering theory is how to detect changes in the constitutive parameters in an inhomogeneous anisotropic medium.
An example are the efforts to detect structural changes in airplane canopies due to prolonged exposure to ultraviolet radiation.
Currently this issue is resolved by simply discarding canopies every few months and replacing them by new ones \cite{CogarColtonMonk:19}.
The difficulty to construct non-destructive evaluation methods for such applications is that
for anisotropic media the corresponding inverse scattering problem does not have a unique solution \cite{GylysC:96}.
A possible remedy to this problem is through the use of qualitative methods \cite{CakoniColtonMonk:11,CakoniColton:14}.
One choice of qualitative properties are so-called target signatures, such as eigenvalues associated with the direct scattering problem that can be reconstructed from the measurement of the scattering data.
Classical methods use the frequency (or its square) as eigenvalue parameter and among these target signatures are scattering frequencies (resonances) and classical transmission eigenvalues.
More recently developed methods use eigenvalue problems with an artificial nonphysical eigenvalue parameter and we refer to \cite{CakoniColtonHaddar:16,CogarColtonMonk:19} for an overview on this subject.
Classical methods which use the frequency as the eigenvalue parameter have two major drawbacks.
Firstly, the measured data is only known for real valued frequencies and thus the
determination of the location of complex eigenvalues is difficult.
However all resonances have a non-zero imaginary part.
The same holds also for classical transmission eigenvalues for absorbing media \cite{CogarColtonMonk:19}.
Secondly, these methods require multi-frequency data, and their measurement-frequencies cannot be chosen arbitrary since they must include the eigenvalues to be measured.
Recently new methods have been reported which overcome both of these drawbacks \cite{CakoniColtonMengMonk:16,CogarColtonMengMonk:17,AudibertCakoniHaddar:17}.
These methods are based on a modified far field operator and  use an artificial nonphysical eigenvalue parameter instead of the frequency.
Hence complex eigenvalues do not give rise to any difficulties and absorbing materials can be treated effectively.
Further, this way the scattering data is required only for a single frequency, and moreover this frequency can be chosen a priori.
Thus these new methods constitute a great practical advance in qualitative inverse scattering methods.
A first variant led to eigenvalue problems of Steklov type \cite{CakoniColtonMengMonk:16}.
In order to increase the sensitivity of the eigenvalues to changes in the material a second variant was proposed \cite{CogarColtonMengMonk:17} through the introduction of an additional sensitivity parameter.
A method similar to \cite{CogarColtonMengMonk:17}, but with a negative index material as background was put forward in \cite{AudibertCakoniHaddar:17}.
This second kind of method leads to equations of transmission type.
Further modifications for the far field operator have been introduced in \cite{CamanoLacknerMonk:17,CogarMonk:20} and \cite{Cogar:20,Cogar:21,Cogar:22} to obtain eigenvalue problems, which are more easily to analyze.
Indeed the analysis of classical transmission, modified transmission and Steklov eigenvalue problems requires sophisticated techniques and can be challenging even for nonabsorbing media~\cite{CakGintHad:10,Halla:21StekAppr,Halla:19StekloffExist,CakoniNguyen:21}.
An understanding of the properties of the eigenvalue problems is essential for the related nondestructive testing methods.
Among the most important questions are the discreteness and existence of eigenvalues, and results on the relation between changes in the material and the corresponding shifting of the eigenvalues.
Clearly  the existence of eigenvalues is essential, because if no eigenvalues exist, then related methods for the inverse scattering problem are useless.
In general for absorbing media the related eigenvalue problems are nonselfadjoint, and the proof of existence of eigenvalues for nonselfadjoint problems is a challenging task.
If the domain and all material coefficients are $C^\infty$-regular, then pseudo-differential calculus with a parameter~\cite{AgranovichVishik:64,Agmon:10} can be used to obtain such results, see e.g.\ \cite{BlastenPaeivaerinta:13,Robbiano:13,CakoniColtonMengMonk:16,CogarColtonMengMonk:17,HaddarMeng:18}.
However, for realistic applications smooth material coefficients cannot be expected.
If the chosen domain of observation and the object under investigation do not coincide, then the material coefficients will be discontinuous at the boundary of the object.
As mentioned in the conclusion of \cite{CakoniColtonMengMonk:16}, such a configuration is of significant practical advantage because computations for the artificial comparison problem can be done ahead of testing.
Also the object might be made of different materials.
Even if the object is homogeneous and the chosen domain and the investigated object coincide, an experienced flaw in the object is unlikely to be totally smooth.
A different technique which works also for rough coefficients/domains is to apply Lidki's Theorem (see e.g.\ \cite[Chapter V, Theorem 6.2]{GohbergKrein:69} for a generalized version).
However, Lidki's Theorem requires the operator to be in a certain Schatten class, which often prevents its application.
For example \cite{HitrikKrupchykOlaPaeivaerinta:11} reports existence results for classical transmission eigenvalues, but does not cover the most interesting case of second order problems in $\setR^3$.
As a remedy, the approach in \cite{Cogar:20,Cogar:21} for Steklov problems was to alter the definition of the particular method such that the corresponding operator is of trace class.
However, such a modification is only known for Steklov problems, but not for transmission problems.
We will see, such a modification for Steklov problems is indeed not necessary to guarantee the existence of eigenvalues.
A second important question, the stability of eigenvalues with respect to changes in the material, was addressed in \cite{CakoniMoskowRome:15,CakoniMoskowRome:18} for classical transmission and in \cite{Cogar:21,Cogar:22} for smoothed Steklov eigenvalue problems.\\

In this article we follow \cite{CamanoLacknerMonk:17} and consider the modified Maxwell Steklov eigenvalue problem in conductive materials to find $(\lambda,\bu)\in\setC\times H(\curl;\Omega)\setminus\{0\}$ such that
\begin{align*}
\curl\bbmu^{-1}\curl \bu- \omega^2 \tilde\bbespilon \bu -i\omega \bbsigma\bu &= 0 \quad\text{in }\Omega,\\
\nv\times\bbmu^{-1}\curl \bu -\lambda \nabla_\Gamma\Delta_\Gamma^{-1}\div_\Gamma(\nv\times\bu) &= 0 \quad\text{on }\partial\Omega,
\end{align*}
with the permeability $\bbmu$, the permittivity $\tilde\bbepsilon$ and the conductivity $\bbsigma$.
The material coefficients $\bbmu$, $\tilde\bbepsilon$, $\bbsigma$ are real, symmetric matrix valued, $\bbsigma$ is nonnegative and $\bbmu$, $\tilde\bbepsilon$ are uniformly positive definite \cite{Monk:03}.
Henceforth we adopt the notation $\bbespilon:=\tilde\bbespilon+\frac{i}{\omega}\bbsigma$ such that $-\omega^2 \tilde\bbespilon-i\omega\bbsigma=-\omega^2\bbepsilon$.
Since we can find an orthonormal basis of $\setR^3$ which diagonalizes $\tilde\bbespilon(x)$ and $\bbsigma(x)$ simultaneously, it follows that $\bbespilon(x)$ is a normal matrix.
The remainder of this article is structured as follows.
In Section~\ref{sec:framework} we recall a Keldysh Theorem to prove the existence of eigenvalues for nonselfadjoint operators, and introduce tools to show that an operator is of finite order.
In Section~\ref{sec:ExistScalar} we apply the former technique to the scalar Steklov eigenvalue problem.
In Section~\ref{sec:ExistMaxwell} we extend the former analysis to the modified Maxwell Steklov problem.
In Section~\ref{sec:stability} we establish in Theorems~\ref{thm:ev-pert}, \ref{thm:improved-rate}, \ref{thm:topnontrivial} the stability of modified Maxwell Steklov eigenvalues with respect to changes in the material.

\section{Notation and framework}\label{sec:framework}

For a vector space $X$ (of scalar functions) let $\boldX:=X^3$ and $\setX:=X^{3\times3}$, whereas $\setN, \setR, \setC$ still denote the sets of natural, real and complex numbers.
For a manifold $D\subset\setR^3$ and a subspace $X\subset L^2(D)$ let $X_*:=\{u\in X\colon \spl u,1\spr_{L^2(D)}=0\}$.
Recall that for Hilbert spaces $X,Y$ a compact operator $K\in L(X,Y)$ is in the Schatten class $K_p(X,Y)$ of order $p\in(0,+\infty)$, if the sequence of singular values $s_n(K), n\in\setN$ of $K$ is $\ell^p(\setN)$-summable.
Further, a compact operator $K\in L(X,Y)$ is said to be of finite order if there exists $p\in(0,+\infty)$ such that $K\in K_p(X,Y)$.
The following theorem, which is a slight reformulation of \cite[Chapter V, Theorem 8.3]{GohbergKrein:69} and \cite[Theorem 4.2, Corollary 3.3]{Markus:88}, will be central for our analysis on the existence of eigenvalues.
\begin{theorem}[Keldysh]\label{thm:Keldysh}
Let $X$ be an infinite dimensional separable Hilbert space and $L(X)$ be the space of bounded operators from $X$ to $X$.
Let $T\in L(X)$ be compact and $K\in L(X)$ be injective, selfadjoint and compact.
Let $T$ or $K$ be of finite order.
Then the spectrum of the operator function $A(\lambda)=I-T-\lambda K$ consists of an infinite sequence of eigenvalues, each with finite algebraic multiplicity, which do not accumulate in $\setC$ and such that for every $\delta>0$ only finitely many eigenvalues lie outside the sectors $\{z\in\setC\colon |\arg z|<\delta\}$ and $\{z\in\setC\colon |\arg z-\pi|<\delta\}$.
If $K$ is nonnegative, then for every $\delta>0$ only finitely many eigenvalues lie outside the sector $\{z\in\setC\colon |\arg z|<\delta\}$.
Further the closure of the space spanned by the generalized eigenspaces of $A(\cdot)$ is dense in $X$.
\end{theorem}
Note that some authors apply Lidki's theorem to establish the existence of eigenvalues.
However, Lidki's theorem requires the operator to be in a particular Schatten class and hence limiting assumptions on the spatial dimension (as in \cite{HitrikKrupchykOlaPaeivaerinta:11}) or the introduction of additional smoothing operators (as in \cite{Cogar:20,Cogar:21}) are applied.
A major difference of Theorem~\ref{thm:Keldysh} to Lidki's Theorem is that $K$ is not required to be in a particular Schatten class, but only needs to be in an arbitrary Schatten class.
Consequently we can avoid restrictive assumptions.
To prove that an operator is of finite order the following two results will be helpful.
\begin{theorem}\label{thm:SobolevFO}
Let $n\in\setN$ and $\Omega\subset\setR^n$ be a bounded Lipschitz domain.
Let $s_2>s_1\geq0$ and $k>\max\{s_2,n/2\}$.
Then the embedding $E_{H^{s_2},H^{s_1}}$ from $H^{s_2}(\Omega)$ into $H^{s_1}(\Omega)$ is compact and in $K_{2k/(s_2-s_1)}(H^{s_2}(\Omega),H^{s_1}(\Omega))$.
\end{theorem}
\begin{proof}
\begin{sloppypar}
It holds $E_{H^k,L^2} \in K_2(H^k(\Omega),L^2(\Omega))$ due to Maurin's Theorem \cite[Theorem 6.61]{AdamsFournier:03}.
Let $(s_n,\varphi_n)_{n\in\setN}$ be the sequence of the normalized eigenpairs of the operator $(E_{H^k,L^2}^*E_{H^k,L^2})^{1/2}$.
Thus $s_n\in l^2(\setN)$.
We denote real interpolation spaces (see e.g.\ \cite{BrennerScott:08}) as $[H_0,H_1]_{K,s,p}$ and Hilbert space interpolation spaces (see e.g.\ \cite[Appendix B]{Bramble:93}) as $[H_0,H_1]_{HS,s}$.
It holds $H^{s}(\Omega)=[L^2(\Omega),H^k(\Omega)]_{K,s/k,2}$ due to \cite[Theorem 14.2.7]{BrennerScott:08}.
The spaces $[L^2(\Omega),H^k(\Omega)]_{K,s/k,2}$ and $[L^2(\Omega),H^k(\Omega)]_{HS,s/k}$ coincide and the norms are equivalent due to \cite[Theorem B.2]{Bramble:93}.
Hence the space $H^{s}(\Omega)$ consists of all functions $u=\sum_{n\in\setN} c_n  \varphi_n$ with $\sum_{n\in\setN} s_n^{2(1-s/k)} |c_n|^2 <+\infty$.
It follows that the eigenvalues of $(E_{H^{s_2},H^{s_1}}^* E_{H^{s_2},H^{s_1}})^{1/2}$ are $\lambda_n=s_n^{(s_2-s_1)/k}$ (whereby $H^s(\Omega)$ is equipped with the scalar product $\spl \cdot,\cdot\spr_{[L^2(\Omega),H^k(\Omega)]_{HS,s/k}}$, $s=s_1,s_2$).
Hence $E_{H^{s_2},H^{s_1}} \in K_{2k/(s_2-s_1)}(H^{s_2}(\Omega),H^{s_1}(\Omega))$.
\end{sloppypar}
\end{proof}
The next lemma is hardly a new result.
However, we only found references for the case $W_1=W_2=W_3$.
\begin{lemma}\label{lem:SCproduct}
\begin{sloppypar}
Let $W_1, W_2$ and $W_3$ be separable Hilbert spaces.
Let $B\in L(W_1,W_2)$ and $A\in L(W_2,W_3)$.
If $s>0$ and either $B\in\SC_\indSC(W_1,W_2)$ or $A\in$ $\SC_\indSC(W_2,W_3)$, then $AB\in\SC_\indSC(W_1,W_3)$.
and
\begin{align*}
\|AB\|_{\SC_\indSC(W_1,W_3)} &\leq \|A\|_{L(W_1,W_2)} \|B\|_{\SC_\indSC(W_1,W_2)}
\quad\text{or}\\
\|AB\|_{\SC_\indSC(W_1,W_3)} &\leq \|A\|_{\SC_\indSC(W_1,W_2)} \|B\|_{L(W_1,W_2)}
\end{align*}
respectively.
\end{sloppypar}
\end{lemma}
\begin{proof}
Due to \cite[Chapter VI Corollary 1.2]{GohbergGoldbergKaashoek:90I} it holds
$s_n(A)=s_n(A^*)$, $s_n(B)=s_n(B^*)$ and $s_n(AB)=s_n(B^*A^*)$.
Thus it suffices to consider the case $B\in\SC_\indSC(W_1,W_2)$.
However, by means of the min-max characterization of singular values it can easily be obtained that $s_n(AB)\leq \|A\|_{L(W_1,W_2)} s_n(B)$ (see also \cite[Chapter VI Proposition 1.3]{GohbergGoldbergKaashoek:90I}) and the claim follows.
\end{proof}

\section{Existence of scalar Steklov eigenvalues}\label{sec:ExistScalar}

In this section we consider the scalar Steklov eigenvalue problem \cite{CakoniColtonMengMonk:16,AudibertCakoniHaddar:17,Harris:21} to find $(\lambda,u)\in\setC\times H^1(\Omega)\setminus\{0\}$ such that
\begin{subequations}\label{eq:Steklov-scalar}
\begin{align}
-\div \bbmu^{-1} \nabla u- \omega^2 \epsilon u &= 0 \quad\text{in }\Omega,\\
\nv\cdot\bbmu^{-1}\nabla u -\lambda u &= 0 \quad\text{on }\partial\Omega.
\end{align}
\end{subequations}
For real $\epsilon\geq0$ and real, positive definite $\bbmu$ the existence of eigenvalues can be studied by means of the spectral theorem for compact, selfadjoint operators.
However, in the nonselfadjoint case existing results~\cite{CakoniColtonMengMonk:16} on the existence of eigenvalues are based on Agmon theory and require the unrealistic assumption that the domain and all parameters are $C^\infty$-smooth.
An alternative approach was proposed in \cite{Cogar:20} to introduce a smoothing operator of sufficiently high degree in the boundary condition of \eqref{eq:Steklov-scalar} and to prove the existence of eigenvalues by means of Lidki's theorem.
However, this does not answer the question if such an approach is indeed necessary, and in addition the introduction of a smoothing operator into the eigenvalue problem seems challenging for other kinds of eigenvalue problems.
\begin{assumption}\label{ass:mu}
Let $\bbmu^{-1} \in \setL^\infty(\Omega)$ be real, symmetric and such that
\begin{align*}
\inf_{x\in\Omega, \xi\in\setC^3, |\xi|=1} \xi^* \bbmu^{-1}(x) \xi >0.
\end{align*}
\end{assumption}
\begin{assumption}\label{ass:NoDirEV}
Let the equation
\begin{align}\label{eq:NoDirichletEV}
-\div \bbmu^{-1} \nabla u- \omega^2 \epsilon u &= 0 \quad\text{in }\Omega,
\qquad
u = 0 \quad\text{on }\partial\Omega,
\end{align}
admit in $H^1(\Omega)$ only the trivial solution $u=0$.
\end{assumption}
\begin{theorem}\label{thm:ExistScalar}
Let $\Omega\subset\setR^3$ be a bounded Lipschitz domain.
Let $\bbmu$ satisfy Assumption~\ref{ass:mu},
$\epsilon \in L^\infty(\Omega)$,
$\omega\in\setR$ and Assumption~\ref{ass:NoDirEV} be satisfied.
Then the spectrum of \eqref{eq:Steklov-scalar} consists of an infinite sequence of eigenvalues, each with finite algebraic multiplicity, which do not accumulate in $\setC$ and such that for every $\delta>0$ only finitely many eigenvalues lie outside the sector $\{z\in\setC\colon |\arg z|<\delta\}$.
\end{theorem}
\begin{proof}
We aim to apply Theorem~\ref{thm:Keldysh}.
To this end we consider $H^1(\Omega)$ with the modified scalar product $\spl u,u' \spr_{H^1(\Omega)}:=\spl \bbmu^{-1}\nabla u,\nabla u' \spr_{\boldL^2(\Omega)}+\spl u,u' \spr_{L^2(\Omega)}$, $u,u'\in H^1(\Omega)$.
Let $T, K\in L(H^1(\Omega))$ be defined through
\begin{align*}
\spl Tu,u'\spr_{H^1(\Omega)}:=\spl (1+\omega^2 \epsilon) u,u' \spr_{L^2(\Omega)}, \quad
\spl Ku,u'\spr_{H^1(\Omega)}:=\spl \tr u,\tr u' \spr_{L^2(\partial\Omega)},
\end{align*}
for all $u,u'\in H^1(\Omega)$ and $A(\lambda):=I-T-\lambda K$.
Then the weak formulation of \eqref{eq:Steklov-scalar} is $\spl A(\lambda)u,u'\spr_{H^1(\Omega)}$ $=0$ for all $u'\in H^1(\Omega)$ or equivalently in operator form $A(\lambda)u=0$.
It readily follows that $K$ is selfadjoint and by means of standard Sobolev embedding theorems that $T$ and $K$ are compact.
With the embedding operator $E_{H^1,H^{3/4}} \in L(H^1(\Omega),H^{3/4}(\Omega))$ and the trace operator $\tr_{H^{3/4},L^2} \in L(H^{3/4}(\Omega),L^2(\partial\Omega))$, $K$ admits the representation $K=E_{H^1,H^{3/4}}^* \tr_{H^{3/4},L^2}^* \tr_{H^{3/4},L^2} E_{H^1,H^{3/4}}$.
Thus it follows with Theorem \ref{thm:SobolevFO} and Lemma \ref{lem:SCproduct} that $K$ is of finite order, though the dimension of $\ker K=H^1_0(\Omega)$ is infinite.
As a remedy we introduce a topological decomposition $H^1(\Omega)=X\oplus^\calT H^1_0(\Omega)$ such that $A(\lambda)$ becomes a \emph{triangular} block operator with respect to this decomposition.
Hence for $u\in H^1(\Omega)$ we define $P_2u\in H^1_0(\Omega)$ to be the solution to
$\spl (I-T) (P_2u-u), u_0'\spr_{H^1(\Omega)}=0$ for all $u_0'\in H^1_0(\Omega)$.
Note that due to Assumption~\ref{ass:NoDirEV} and the Fredholm alternative the operator $P_2$ is indeed well-defined and bounded $P_2\in L(H^1(\Omega),H^1_0(\Omega))$.
It is also straightforward to see that $P_2$ is a projection.
Consequently with $X:=\ran P_1$, $P_1:=(I-P_2)$ we constructed a topological decomposition $H^1(\Omega)=X\oplus^\tau H^1_0(\Omega)$ with associated projections $P_1, P_2$ such that $A(\lambda)$ admits a block operator representation
\begin{align*}
A(\lambda) \approx
\bpm P_1^o A(\lambda) |_X & P_1^o A(\lambda) |_{H^1_0(\Omega)} \\ 0 & P_2^o A(\lambda) |_{H^1_0(\Omega)} \epm,
\end{align*}
whereat $P_1^0, P_2^o$ are the orthogonal projections onto $X, H^1_0(\Omega)$.
Due to Assumption~\ref{ass:NoDirEV} and the Fredholm alternative the operator $P_2^o A(\lambda) |_{H^1_0(\Omega)}$ is bijective.
Hence the eigenvalues of $A(\lambda)$ and $A_X(\lambda):=P_1^o A(\lambda) |_X=I_X-P_1^oT|_X-\lambda P_1^oK|_X$ coincide.
Note that $\{\tr u\colon u\in X\}=\{\tr u\colon u\in H^1(\Omega)\}=H^{1/2}(\partial\Omega)$ from which we deduce that $X$ has an infinite dimension.
Now $P_1^oK|_X$ $=$ $P_1^oKP_1|_X$ is injective and Theorem~\ref{thm:Keldysh} can be applied.
Hence the proof is finished.
\end{proof}
We close this section with a discussion on Assumption~\ref{ass:NoDirEV}.
\begin{remark}[Validity of Assumption~\ref{ass:NoDirEV}]\label{rem:assinj-val}
Let us discuss the validity of the Assumption~\ref{ass:NoDirEV}.
We are interested mainly in nonreal $\epsilon$.
In this case $\epsilon$ is typically of the form $\epsilon=\tilde\epsilon+\frac{i}{\omega}\sigma$ with uniformly positive $\tilde\epsilon\in L^\infty(\Omega)$ and $\sigma\in L^\infty(\Omega)$ with $\sigma\geq0$ in $\Omega$.
If $\sigma$ vanishes nowhere, then the injectivity of \eqref{eq:NoDirichletEV} follows by elementary means.
If only $\sigma>0$ in a subset $\Omega_1\subset\Omega$, then it only follows that any solution $u\in H^1(\Omega)$ to \eqref{eq:NoDirichletEV} satisfies $u=0$ in $\Omega_1$.
Thence if $\bbmu^{-1}, \epsilon, \Omega$ allow the application of a unique continuation principle, then it follows $u=0$ in $\Omega$ as well.
In cases which are not covered by the aforementioned considerations one can argue that there exist at most countably many frequencies $\omega$ for which the Assumption~\ref{ass:NoDirEV} is violated.
\end{remark}
\begin{remark}[Necessity of Assumption~\ref{ass:NoDirEV}]\label{rem:assinj-nec}
When comparing the assumptions of Theorem~\ref{thm:ExistScalar} with the assumptions of similar works (e.g.\ \cite{CakoniColtonMengMonk:16,AudibertCakoniHaddar:17}), the question appears if Assumption~\ref{ass:NoDirEV} is indeed necessary.
We do not have a definite answer to this question.
However, we note that when Steklov eigenvalues are studied by means of a Neumann-to-Dirichlet or a Robin-to-Dirichlet operator $R\in L^2(\partial\Omega)$, the Assumption~\ref{ass:NoDirEV} is implicitly used (often without mentioning it).
Indeed for an eigenpair $(\lambda,u)$ of \eqref{eq:Steklov-scalar} it is derived that $\lambda R\tr u=\tr u$.
To deduce for $\lambda\neq0$ that $\lambda^{-1}$ is an eigenvalue of $R$ it is then necessary that $\tr u\neq0$.
We do not know how $\tr u\neq0$ can be derived without any further assumptions.
However, if Assumption~\ref{ass:NoDirEV} is satisfied, then a solution $(\lambda,u)$ of \eqref{eq:Steklov-scalar} with $\tr u=0$ would contradict $u\neq0$.
On the other hand to deduce the existence of infinitely man distinct eigenvalues of \eqref{eq:Steklov-scalar} from an analysis of $R$ the Assumption~\ref{ass:NoDirEV} is also necessary.
For example if $\epsilon$ is real, then $R$ is selfadjoint and the existence of an orthonormal eigenbasis of $R$ follows.
Since $L^2(\partial\Omega)$ is infinite dimensional, there exists an infinite sequence of eigenvalues $\tau_n, n\in\setN$ to $R$.
However, in the case that the codimension of $\ker R$ is finite only finite many eigenvalues $\tau_n,n=1,\dots,N\in\setN$ are nonzero.
Hence, we can only deduce a finite dimensional sequence of eigenvalues $\tau_n^{-1},n=1,\dots,N$ to \eqref{eq:Steklov-scalar}.
However, if Assumption~\ref{ass:NoDirEV} is satisfied, then it follows that $R$ is injective and the former situation cannot occur.
See e.g.~\cite{AudibertCakoniHaddar:17,Cogar:20,Cogar:22} for works which impose an assumption identical or similar to Assumption~\ref{ass:NoDirEV}.
\end{remark}
\begin{remark}[Similar assumptions to Assumption~\ref{ass:NoDirEV}]\label{rem:assinj-sim}
We mention that in some works \cite{CakoniColtonMengMonk:16,Harris:21} on scalar Steklov eigenvalues the assumption that the equation
\begin{align*}
-\div \bbmu^{-1} \nabla u- \omega^2 \epsilon u &= 0 \quad\text{in }\Omega,
\qquad
\nv\cdot\bbmu^{-1}\nabla u = 0 \quad\text{on }\partial\Omega,
\end{align*}
admits only the trivial solution $u=0$, appears (i.e.\ the Dirichlet boundary condition in \eqref{eq:NoDirichletEV} is replaced by a Neumann boundary condition).
This assumption is used to introduce a Neumann-to-Dirichlet operator, and subsequently to transform the Steklov eigenvalue problem to an eigenvalue problem posed on the space $L^2(\partial\Omega)$.
Note that instead one can work with a Robin-to-Dirichlet operator (with a suitably chosen Robin parameter) as in e.g.\ in \cite{Cogar:20} to avoid such an assumption.
\end{remark}

\section{Existence of modified Maxwell Steklov eigenvalues}\label{sec:ExistMaxwell}

In this section we consider the modified Maxwell Steklov eigenvalue problem \cite{CamanoLacknerMonk:17,Halla:21StekAppr,GongSunWU:21,Gong:22} to find $(\lambda,\bu)\in\setC\times H(\curl;\Omega)\setminus\{0\}$ such that
\begin{subequations}\label{eq:Steklov-em}
\begin{align}
\curl\bbmu^{-1}\curl \bu- \omega^2 \bbespilon \bu &= 0 \quad\text{in }\Omega,\\
\nv\times\bbmu^{-1}\curl \bu +\lambda\, S\bu &= 0 \quad\text{on }\partial\Omega.
\end{align}
\end{subequations}
Here the operator $S$ which appears in the boundary condition is as follows.
We denote the differential operators on $\partial\Omega$ with a subscript $\Gamma$.
Thence
\begin{align*}
S:=-\nabla_\Gamma\Delta_\Gamma^{-1}\div_\Gamma\tr_{\nv\times},
\end{align*}
whereat $\tr_{\nv\times}\bu:=\nv\times \tr \bu$ and $\Delta_\Gamma^{-1}$ maps into $H^1_*(\partial\Omega)$.
Note that due to
\begin{align*}
\tr_{\nv\times}\in L(H(\curl;\Omega), H^{-1/2}(\div_\Gamma;\partial\Omega))
\end{align*}
the operator $S$ is indeed bounded: $S\in L(H(\curl;\Omega), \boldL^2(\partial\Omega))$.
The eigenvalue problem~\eqref{eq:Steklov-em} was introduced in \cite{CamanoLacknerMonk:17} as a modification of the standard Steklov eigenvalue problem ($S=I$) to cope with the complicated nature of this equation.
The Fredholmness of the standard Steklov eigenvalue problem was later established in \cite{Halla:21StekAppr} and the distribution of eigenvalues in the selfadjoint case reported in \cite{Halla:19StekloffExist}.
The approach to introduce a smoothing was generalized from the scalar case \cite{Cogar:20} in \cite{Cogar:21} to \eqref{eq:Steklov-em}.
Note that our forthcoming analysis can easily be extended to those problems.
We are going to establish the existence of eigenvalues of \eqref{eq:Steklov-em} in Theorem~\ref{thm:ExistMaxwell}.
Before that we recall some embedding results for the spaces
\begin{align*}
H(\curl,\div\bbespilon,\tr_{\nv\times}^0;\Omega)&:=\{\bu\in H(\curl;\Omega)\colon \div(\bbespilon\bu)\in L^2(\Omega), \, \nv\times\bu=0 \text{ on }\partial\Omega\},\\
H(\curl,\div\bbespilon,\tr_{\nv\cdot\bbespilon}^0;\Omega)&:=\{\bu\in H(\curl;\Omega)\colon \div(\bbespilon\bu)\in L^2(\Omega), \, \nv\cdot \bbespilon\bu=0 \text{ on }\partial\Omega\},
\end{align*}
which are equipped with the norm
\begin{align*}
\|\bu\|^2_{H(\curl,\div\bbespilon;\Omega)}:=\|\bu\|^2_{\boldL^2(\Omega)}+\|\curl\bu\|^2_{\boldL^2(\Omega)}+\|\div(\bbespilon\bu)\|^2_{L^2(\Omega)},
\end{align*}
and investigate the properties of the operator $S$.
\begin{assumption}\label{ass:epscoercive}
Let $\bbespilon \in \setL^\infty(\Omega)$ be such that
$\bbespilon_-:=
\inf_{x\in\Omega, \xi\in\setC^3, |\xi|=1} \Re(\xi^* \bbespilon^{-1}(x) \xi) >0$.
\end{assumption}
\begin{proposition}\label{prop:Weber}
Let $\bbespilon$ satisfy Assumption~\ref{ass:epscoercive}.
Then the spaces $H(\curl,\div\bbespilon,\tr_{\nv\times}^0;\Omega)$ and $H(\curl,\div\bbespilon,\tr_{\nv\cdot\bbespilon}^0;\Omega)$ embed compactly into $\boldL^2(\Omega)$.
\end{proposition}
\begin{proof}
For real valued matrix functions $\bbespilon$ this result is well known due to \cite{Weber:80}.
The proof from \cite{Weber:80} for complex valued matrix functions $\bbespilon$ which satisfy Assumption~\ref{ass:epscoercive} requires only one minor adaptation: all inequalities of the kind $\bbespilon_- \|\bu\|^2_{\boldL^2(\Omega)} \leq \spl \bbespilon \bu,\bu \spr_{\boldL^2(\Omega)}$ have to  be replaced by $\bbespilon_- \|\bu\|^2_{\boldL^2(\Omega)} \leq \Re(\spl \bbespilon \bu,\bu \spr_{\boldL^2(\Omega)})$.
See also the proofs of \cite[Theorems~8.1.1,~8.1.3]{AssousCiarletLabrunie:18}.
\end{proof}
The set $\calP:=\{\Omega_j\}, j=1,...,N$ is called a partition of $\Omega$ if $\Omega_j\subset\Omega, j=1,...,N$ are mutually disjoint connected Lipschitz domains, and it holds that $\ol{\Omega}=\bigcup_{j=1}^N \ol{\Omega_j}$.
Let
\begin{align*}
PW^{1,\infty}(\Omega,\calP):=\{\epsilon\in L^\infty(\Omega)\colon \epsilon|_{\Omega_j}\in W^{1,\infty}(\Omega_j), j=1,\dots,N\}.
\end{align*}
We make an assumption on the topology of $\Omega$ as in \cite{Ciarlet:20} and refer to \cite[Chapter~3, Remark 3.2.1, Theorem~3.3.2]{AssousCiarletLabrunie:18} for a discussion on this assumption.
\begin{assumption}\label{ass:top}
The Lipschitz domain $\Omega$ satisfies one of the following topology conditions:
\begin{itemize}
 \item $\mathbf{(Top)}_{I=0}\colon$ ``given any curl-free vector field $\bu\in\boldC^1(\Omega)$, there exists $p\in C(\Omega)$ such that $\bu=\nabla p$ in $\Omega$'';
 \item $\mathbf{(Top)}_{I>0}\colon$ ``there exist $I\in\setN$ nonintersecting manifolds, $\Sigma_1,\dots,\Sigma_I$, with boundaries $\partial\Sigma_i\subset\partial\Omega$, such that, if we let $\dot\Omega:=\Omega\setminus\bigcup_{i=1}^I\Sigma_i$, given any curl-free vector field $\bu\in C(\Omega)$, there exists $\dot p\in C(\dot\Omega)$ such that $\bu=\nabla\dot p$ in $\dot\Omega$”.
\end{itemize}
\end{assumption}
We call $\bbespilon\in\setL^{\infty}(\Omega)$ normal, if $\bbespilon(x)$ is a normal matrix for almost all $x\in\Omega$.
\begin{proposition}\label{prop:Ciarlet}
\begin{sloppypar}
Let $\Omega$ be a Lipschitz domain which satisfies the topology Assumption \ref{ass:top}.
Let $\bbespilon\in\setP\setW^{1,\infty}(\Omega,\calP)$ be normal and satisfy Assumption~\ref{ass:epscoercive}.
Then there exist constants $s_1, s_2>0$ such that $H(\curl,\div\bbespilon,\tr_{\nv\times}^0;\Omega)$ and $H(\curl,\div\bbespilon,\tr_{\nv\cdot\bbespilon}^0;\Omega)$ embed continuously in $\boldH^s(\Omega)$ for all $s\in(0,s_1)$ and $s\in(0,s_2)$ respectively.
\end{sloppypar}
\end{proposition}
\begin{proof}
The result follows from \cite{Ciarlet:20}.
Let $\bu\in H(\curl,\div\bbespilon,\tr_{\nv\times}^0;\Omega)$ or $\bu\in H(\curl,\div\bbespilon,\tr_{\nv\cdot\bbespilon}^0;\Omega)$.
Due to Theorems~5.5 and~5.6 of \cite{Ciarlet:20} we can make decomposition $\bu=\boldv_\mathrm{reg}+\nabla p$ such that $\|\boldv_\mathrm{reg}\|_{\boldH^{1/2}(\Omega)} \leq C \|\bu\|_{H(\curl,\div\bbespilon;\Omega)}$ with a constant $C>0$ independent of $\bu$.
The function $p\in H^1(\Omega)$ is a solution to (5.4) or (5.7) in \cite{Ciarlet:20} respectively.
Due to Lemma~6.10 of \cite{Ciarlet:20} there exists $s>0$ such that the $H^{s-1}(\Omega)$-norm of the right hand-sides in (5.4) and (5.7) is bounded by a constant times $\|\bu\|_{H(\curl,\div\bbespilon;\Omega)}$.
Hence Theorem~6.8 of \cite{Ciarlet:20} yields that $\|p\|_{H^{1+s}(\Omega)}$ is bounded by a constant times $\|\bu\|_{H(\curl,\div\bbespilon;\Omega)}$.
Note that Lemma~6.10 and Theorem~6.8 of \cite{Ciarlet:20} are formulated only for matrix functions $\bbespilon$ which are multiples of the identity matrix.
However, as the explained on Page~3029 of \cite{Ciarlet:20} Lemma~6.10 and Theorem~6.8 of \cite{Ciarlet:20} extend to normal matrix functions $\bbespilon$.
Thus the proof is finished.
\end{proof}
Next we investigate the properties of the operator $S$.
\begin{lemma}\label{lem:Scompact}
The operator $S\in L(H(\curl;\Omega),\boldL^2(\partial\Omega))$ is compact.
\end{lemma}
\begin{proof}
We recall $S=-\nabla_\Gamma\Delta_\Gamma^{-1}\div_\Gamma\tr_{\nv\times}$ and that $-\Delta_\Gamma^{-1}f$ is the solution to the problem to find $u\in H^1_*(\partial\Omega)$ such that
\begin{align*}
\spl\nabla_\Gamma u,\nabla_\Gamma u'\spr_{\boldL^2(\partial\Omega)}
=\spl f,u'\spr_{H^{-1}(\partial\Omega)\times H^{1}(\partial\Omega)}
\quad\text{for all}\quad u'\in H^1_*(\partial\Omega).
\end{align*}
For $\bu\in H(\curl;\Omega)$ we know that $\div\tr_{\nv\times}\bu\in H^{-1/2}(\partial\Omega)$, and hence one might expect that $S\bu$ admits some extra regularity, i.e.\ $S\bu\in H^{1+s}(\partial\Omega)$ with $s>0$.
However, for Lipschitz domains the usual definition of $H^{s}(\partial\Omega)$ via local coordinates is only well-defined for $s\in[-1,1]$.
Nonstandard definitions of $H^{s}(\partial\Omega)$, $s>1$ e.g.\ as $\tr H^{s+1/2}(\Omega)$ exist, but such spaces lack some classical properties of Sobolev spaces.
In the following we consider the space $H^1_*(\partial\Omega)$ equipped with the norm $\|\nabla_\Gamma u\|_{\boldL^2(\partial\Omega)}$ (which is equivalent to the standard  $H^1(\partial\Omega)$-norm on $H^1_*(\partial\Omega)$).
Will work with the spaces $[L^2_*(\partial\Omega),H^1_*(\partial\Omega)]_{HS,s}$, $s\in\setR$ generated by Hilbert space interpolation (see e.g.\ \cite[Appendix B]{Bramble:93}).
For $s\in[-1,1]$ it holds $[L^2_*(\partial\Omega),H^1_*(\partial\Omega)]_{HS,s}=H^s_*(\partial\Omega)$ (see e.g.\ \cite[p.\ 232]{Agranovich:15}).
However, the spaces $[L^2_*(\partial\Omega),H^1_*(\partial\Omega)]_{HS,s}$ are also well-defined for $s>1$ (and in this case the name \emph{extra}polation spaces might be more apt).
Let $(s_n,\varphi_n), n\in\setN$ be the eigenpairs of $(E^*E)^{1/2}$ with the embedding operator $E\in L(H^1_*(\partial\Omega),L^2_*(\partial\Omega))$.
Recall that $[L^2_*(\partial\Omega),H^1_*(\partial\Omega)]_{HS,s}$ consists of all functions $u=\sum_{n\in\setN} c_n \varphi_n$ with $\sum_{n\in\setN} s_n^{2(1-s)} |c_n|^2 <+\infty$.
Since $\div\tr_{\nv\times}\bu\in H^{-1/2}(\partial\Omega)$ it follows $\div\tr_{\nv\times}\bu \in [L^2_*(\partial\Omega),H^1_*(\partial\Omega)]_{HS,-1/2}$ and thus $-\Delta_\Gamma^{-1}\div\tr_{\nv\times}\bu\in [L^2_*(\partial\Omega),H^1_*(\partial\Omega)]_{HS,3/2}$.
Thus $S\bu=\sum_{n\in\setN} c_n \nabla_\Gamma\varphi_n$ with $\sum_{n\in\setN} s_n^{-1} |c_n|^2 <+\infty$.
Since $\|S\bu\|_{\boldL^2(\partial\Omega)}^2=\sum_{n\in\setN} |c_n|^2$ it follows that $S\in L(H(\curl;\Omega),\boldL^2(\partial\Omega))$ is compact.
\end{proof}
\begin{lemma}\label{lem:S-Hs}
Let $\Omega$ be a $C^\infty$ domain.
Then it holds that
\begin{align*}
S\in L(H(\curl;\Omega),\boldH^{1/2}(\partial\Omega)).
\end{align*}
\end{lemma}
\begin{proof}
For $C^\infty$ domains the Sobolev spaces $H^s(\partial\Omega)$, $s\in\setR$ can be defined in the convenient way and it follows from \cite[Proposition~1]{delaBourdonnaye:93} that indeed $S\in L(H(\curl;\Omega),\boldH^{1/2}(\partial\Omega))$.
\end{proof}
\begin{assumption}\label{ass:NoDirEV-MW}
For $\ker S:=\{\bu\in H(\curl;\Omega)\colon S\bu=0\}$ let the equation
\begin{align*}
\spl \bbmu^{-1}\curl\bu,\curl\bu'\spr_{\boldL^2(\Omega)}
-\omega^2 \spl \bbespilon\bu,\bu'\spr_{\boldL^2(\Omega)} =0
\quad\text{for all }\bu'\in\ker S
\end{align*}
admit only the trivial solution $\bu=0$ in $\ker S$.
\end{assumption}
Now we are prepared to establish the main result of this section.
\begin{theorem}\label{thm:ExistMaxwell}
Let $\Omega\subset\setR^3$ be a bounded Lipschitz domain.
Let $\bbmu$ satisfy Assumption~\ref{ass:mu} and $\bbespilon$ satisfy Assumption~\ref{ass:epscoercive}.
Let $\Omega$ be a $C^\infty$ domain, or $\Omega$ satisfy the topology Assumption~\ref{ass:top} and $\bbespilon\in\setP\setW^{1,\infty}(\Omega)$ be normal.
Let $\omega\in\setR\setminus\{0\}$ and Assumption~\ref{ass:NoDirEV-MW} be satisfied.
Then the spectrum of \eqref{eq:Steklov-em} consists of an infinite sequence of eigenvalues, each with finite algebraic multiplicity, which do not accumulate in $\setC$ and such that for every $\delta>0$ only finitely many eigenvalues lie outside the sector $\{z\in\setC\colon |\arg z|<\delta\}$.
\end{theorem}
\begin{proof}
\begin{sloppypar}
We aim to mimic the technique for the scalar case (see Theorem~\ref{thm:ExistScalar}) and to apply Theorem~\ref{thm:Keldysh}.
To this end we consider $H(\curl;\Omega)$ equipped with the modified scalar product $\spl \bu,\bu' \spr_{H(\curl;\Omega)}:=\spl \bbmu^{-1}\curl \bu,\curl \bu' \spr_{\boldL^2(\Omega)}$ $+$ $\spl \bu,\bu' \spr_{\boldL^2(\Omega)}$, $\bu,\bu'\in H(\curl;\Omega)$.
Let $T, K\in L(H(\curl;\Omega))$ be defined through
\begin{align*}
\spl T\bu,\bu'\spr_{H(\curl;\Omega)}:=\spl (1+\omega^2 \bbespilon) \bu,\bu' \spr_{\boldL^2(\Omega)}, \hspace{1.2mm}
\spl K\bu,\bu'\spr_{H(\curl;\Omega)}:=\spl S\bu,S\bu' \spr_{\boldL^2(\partial\Omega)},
\end{align*}
\end{sloppypar}
\noindent
for all $\bu,\bu'\in H(\curl;\Omega)$ and $A(\lambda):=I-T-\lambda K$.
Then the weak formulation of \eqref{eq:Steklov-em} is $\spl A(\lambda)\bu,\bu'\spr_{H(\curl;\Omega)}=0$ for all $\bu'\in H(\curl;\Omega)$ or equivalently in operator form $A(\lambda)\bu=0$.
Next we introduce a suitable topological decomposition $H(\curl;\Omega)=X\oplus^\calT \ker S$ (with $\ker S\subset H(\curl;\Omega)$) such that $A(\lambda)$ becomes a \emph{triangular} block operator with respect to this decomposition.
For $\bu\in H(\curl;\Omega)$ we define $P_2\bu\in \ker S$ to be the solution to $\spl (I-T) (P_2\bu-\bu), \bu'\spr_{H(\curl;\Omega)}=0$ for all $\bu'\in \ker S$, or equivalently formulated $P_2\bu\in\ker S$ solves
\begin{align}\label{eq:P2}
\spl \bbmu^{-1}\curl (P_2\bu-\bu),\curl \bu'\spr_{\boldL^2(\Omega)}-\omega^2\spl \bbespilon (P_2\bu-\bu),\bu'\spr_{\boldL^2(\Omega)}=0
\end{align}
for all $\bu'\in\ker S$.
Note that $Z:=\{\nabla z\colon z\in H^1_*(\Omega)\}\subset\ker S$.
Let
\begin{align*}
Y&:=\{\bu\in\ker S\colon \spl \bu,\nabla z\spr_{H(\curl;\Omega)}=0 \text{ for all }z\in H^1_*(\Omega)\}\\
&=\{\bu\in\ker S\colon \div \bu=0\text{ in }\Omega, \, \nv\cdot \bu=0 \text{ on }\partial\Omega\}
\end{align*}
be the orthogonal complement of $Z$ in $\ker S$.
Let $P_Y$ and $P_Z$ the orthogonal projections on the respective subspaces.
Then with $\tilde T:=P_Y-P_Z$ the
sesquilinearform $\spl (I-T) \bu, \bu'\spr_{H(\curl;\Omega)}$, $\bu,\bu'\in \ker S$ is weakly $\tilde T$-coercive.
Hence Assumption~\ref{ass:NoDirEV-MW} and the Fredholm alternative yield that the operator associated to the sesquilinearform $\spl (I-T) \bu, \bu'\spr_{H(\curl;\Omega)}$, $\bu,\bu'\in \ker S$ is bijective.
Thus $P_2\bu$ is indeed well-defined and bounded $P_2\in L(H(\curl;\Omega),\ker S)$.
Further it can easily be seen that $P_2$ is a projection.
Consequently with $X:=\ran P_1$, $P_1:=(I-P_2)$ we constructed a topological decomposition $H(\curl;\Omega)=X\oplus^\calT \ker S$ with associated projections $P_1, P_2$ such that $A(\lambda)$ admits a block operator representation
\begin{align*}
A(\lambda) \approx
\bpm P_1^o A(\lambda) |_X & P_1^o A(\lambda) |_{\ker S} \\ 0 & P_2^o A(\lambda) |_{\ker S} \epm,
\end{align*}
whereat $P_1^o, P_2^o$ are the orthogonal projections onto $X,\ker S$.
As we already discussed the operator $P_2^o A(\lambda) |_{\ker S}$ is bijective.
Hence the eigenvalues of $A(\lambda)$ and $A_X(\lambda):=P_1^o A(\lambda) |_X=I_X-P_1^oT|_X-\lambda P_1^oK|_X$ coincide.
Due to the definition of the space $X$ the operator $P_1^oK|_X$ is injective.
It can easily be seen that $P_1^oK|_X$ is selfadjoint, and $P_1^oK|_X$ is compact due to Lemma~\ref{lem:Scompact}.
Note that it follows from~\eqref{eq:P2} and $Z\subset \ker S$ that $X\subset H(\curl,\div\bbespilon,\tr_{\nv\cdot\bbespilon}^0;\Omega)$.
Thus Proposition~\ref{prop:Weber} implies that $P_1^oT|_X$ is compact.
If $\Omega$ is a $C^\infty$ domain, then $P_1^oK|_X$ is of finite order due to Lemma~\ref{lem:S-Hs}, Theorem \ref{thm:SobolevFO} and Lemma \ref{lem:SCproduct}.
If $\Omega$ satisfies the topology Assumption~\ref{ass:top} and $\bbespilon\in\setP\setW^{1,\infty}(\Omega)$ is normal, then $P_1^oT|_X$ is of finite order due to Proposition~\ref{prop:Ciarlet}, Theorem~\ref{thm:SobolevFO} and Lemma~ \ref{lem:SCproduct}.
It remains to show that $X$ has an infinite dimension.
To this end we apply some results from \cite{BuffaCostabelSheen:02}.
Due to \cite[Theorem~4.1]{BuffaCostabelSheen:02} we know that the image of $H(\curl;\Omega)$ under the tangential trace operator $\tr_{\nv\times}$ equals $H^{-1/2}(\div_\Gamma;\partial\Omega)$.
Let $\calH(\partial\Omega):=\{p\in H^1_*(\partial\Omega)\colon\Delta_\Gamma p\in H^{-1/2}_*(\partial\Omega)\}$.
Due to \cite[Theorem~5.5]{BuffaCostabelSheen:02} it holds that $\nabla_\Gamma\calH(\partial\Omega)\subset H^{-1/2}(\div_\Gamma;\partial\Omega)$.
For $p\in\calH(\partial\Omega)$ let $\bu(p)\in H(\curl;\Omega)$ be such that $\nv\times\bu(p)=\nabla_\Gamma p$ on $\partial\Omega$.
It follows that $P_1\bu(p)\in X$ and that $SP_1\bu(p)=\nabla_\Gamma p$.
Hence $\nabla_\Gamma\calH(\partial\Omega)\subset SX$.
Since the dimension of $\nabla_\Gamma\calH(\partial\Omega)$ is infinite, it follows that the dimension of $X$ is infinite too.
Now Theorem~\ref{thm:Keldysh} can be applied.
Thus the proof is finished.
\end{proof}
To interpret the injectivity Assumption~\ref{ass:NoDirEV-MW} we can apply similar considerations as in Remarks~\ref{rem:assinj-val}, \ref{rem:assinj-nec}, \ref{rem:assinj-sim} for the scalar case.
Compared to the Assumption~\ref{ass:NoDirEV} for the scalar case, Assumption~\ref{ass:NoDirEV-MW} looks a bit odd.
However, this is exactly the Assumption which is e.g.\ missing in \cite{CamanoLacknerMonk:17}.

\section{Stability of modified electromagnetic Steklov eigenvalues}\label{sec:stability}

In this section we consider eigenvalue problem \eqref{eq:Steklov-em} for different pairs $(\bbmu,\bbespilon)=(\bbmu_0,\bbespilon_0), (\bbmu_h,\bbespilon_h)$ of parameters, and investigate the convergence of eigenvalues as $(\bbmu_h,\bbespilon_h)$ $\to$ $(\bbmu_0,\bbespilon_0)$.
That is we study the stability of eigenvalues of \eqref{eq:Steklov-em} with respect to perturbations in $(\bbmu,\bbespilon)$.
If $(\bbmu_h,\bbespilon_h)$ converges to $(\bbmu_0,\bbespilon_0)$ in $\setL^\infty(\Omega)$-norms, then it is rather straightforward to obtain the convergence of eigenvalues.
However, the former consideration excludes an important class of perturbations.
Consider e.g.\ $\bbmu_h=\bbmu_0+\delta_1\chi_{B_h(x_1)}\setI$, $\bbespilon_h=\bbespilon_0+\delta_2\chi_{B_h(x_2)}\setI$ with scalar constants $\delta_1,\delta_2$, $\chi_B$ being the characteristic function of a set $B$, $\setI$ being the $3\times3$ identity matrix, $x_1,x_2\in\Omega$, $h\geq0$ and the balls $B_h(x):=\{y\in\setR^3\colon |x-y|<h\}$.
Then $\|\bbmu_0-\bbmu_h\|_{\setL^\infty(\Omega)}=|\delta_1|$, $\|\bbespilon_0-\bbespilon_h\|_{\setL^\infty(\Omega)}=|\delta_2|$, while $\|\bbmu_0-\bbmu_h\|_{\setL^p(\Omega)}=O(|\delta_1|h^{3/p})$, $\|\bbespilon_0-\bbespilon_h\|_{\setL^p(\Omega)}=O(|\delta_2|h^{3/p})$ for each $p\in[1,\infty)$.
Hence in this case $(\bbmu_h,\bbespilon_h)$ converges to $(\bbmu_0,\bbespilon_0)$ in $\setL^p(\Omega)$, $p\in[1,\infty)$, but not in $\setL^\infty(\Omega)$.
Thus we would like to obtain stability results in terms of $\setL^p(\Omega)$-norms.
To this end we are going to reformulate the eigenvalue problems \eqref{eq:Steklov-em} into eigenvalue problems for analytic Fredholm operator functions $\setT_0(\cdot), \setT_h(\cdot)$.
Subsequently we will estimate in Lemma~\ref{lem:norm-est} the difference $\|\setT_0(\lambda)- \setT_h(\lambda)\|_{L(V)}$ in terms of $\|\bbmu_0-\bbmu_h\|_{\setL^p(\Omega)}$ and $\|\bbespilon_0-\bbespilon_h\|_{\setL^p(\Omega)}$.
Then general perturbation theory for analytic operator functions \cite{Karma:96a,Karma:96b} will yield in Theorem~\ref{thm:ev-pert} the convergence of eigenvalues in terms of $\|\bbmu_0-\bbmu_h\|_{\setL^p(\Omega)}$ and $\|\bbespilon_0-\bbespilon_h\|_{\setL^p(\Omega)}$.
In Theorem~\ref{thm:improved-rate} we will apply a more specialized perturbation theorem \cite{Moskow:15} to obtain an improved convergence rate under an additional assumption.
In Theorem~\ref{thm:topnontrivial} we will extend the former results to topologically nontrivial domains.
Note that the forthcoming analysis which will result in Theorems~\ref{thm:ev-pert}, \ref{thm:improved-rate} and \ref{thm:topnontrivial} is a significant improvement of \cite{Cogar:21,Cogar:22}.
Indeed here we will require only uniform $\setL^\infty$ bounds on $\bbmu_h,\bbespilon_h$ (see Assumption~\ref{ass:mueps}).
In comparison \cite{Cogar:21} requires that $\|\bbespilon_0-\bbespilon_h\|_{\setL^\infty(\Omega)}$ is sufficiently small and the analysis \cite{Cogar:22} for scalar problems requires that $A_h$ converges to $A_0$ in $\setP\setW^{1,\infty}(\Omega,\calP)$ with a fixed partition $\calP$, whereat $A_0,A_h$ in \cite{Cogar:22} correspond to $\bbmu_0^{-1},\bbmu^{-1}_h$ here.
\begin{assumption}\label{ass:mueps}
\begin{sloppypar}
Let $(\bbmu_h)_{h\geq0}, (\bbespilon_h)_{h\geq0}$ be sequences with $\bbmu_h,\bbespilon_h\in\setL^\infty(\Omega)$ such that $\sup_{h\geq0}\|\bbmu_h^{-1}\|_{\setL^\infty(\Omega)}$, $\sup_{h\geq0}\|\bbespilon_h\|_{\setL^\infty(\Omega)}<\infty$.
Let $(\bbmu_h)_{h\geq0}$ be real and symmetric, $\bbespilon_0$ be normal, and
$\mu_-:=\hspace{-0mm}\inf_{x\in\Omega, \xi\in\setC^3, |\xi|=1, h\geq0} \hspace{-0mm} \xi^* \bbmu_h^{-1}(x) \xi >0$
and
$\epsilon_-:=\hspace{-0mm}\inf_{x\in\Omega, \xi\in\setC^3, |\xi|=1, h\geq0} \hspace{-0mm} \Re(\xi^* \bbespilon_h(x) \xi) >0$.
\end{sloppypar}
\end{assumption}
Let
\begin{align*}
V_h&:=H(\curl,\div\bbespilon_h^0,\tr_{\nv\cdot\bbespilon_h}^0;\Omega):=\{\bu\in H(\curl;\Omega)\colon \div(\bbespilon_h\bu)=0 \text{ in }\Omega, \, \nv\cdot\bbespilon_h\bu=0\text{ on }\partial\Omega\},\\
\spl \bu,\bu'\spr_{V_h}&:=\spl\bbmu_h^{-1}\curl\bu,\curl\bu'\spr_{\boldL^2(\Omega)}.
\end{align*}
Let $V:=V_0$.
For $(\bbmu_h)_{h\geq0}, (\bbespilon_h)_{h\geq0}$ satisfying Assumption~\ref{ass:mueps} it is well-known that we can make the topological decomposition $H(\curl;\Omega)=V_h\oplus^\calT \nabla H^1_*(\Omega)$.
The projection onto $V_h$ is given by
\begin{align*}
\PVh\bu:=\bu-\nabla w
\end{align*}
with $w\in H^1_*(\Omega)$ being the solution to
\begin{align}\label{eq:eqforw}
\spl \bbespilon_h \nabla w, \nabla w' \spr_{\boldL^2(\Omega)} = \spl \bbespilon_h \bu, \nabla w' \spr_{\boldL^2(\Omega)} \quad\text{for all } w\in H^1_*(\Omega).
\end{align}
In addition if $\Omega$ is simply connected, which we assume from now on, then $\spl\cdot,\cdot\spr_{V_h}$ is on $V_h$ an equivalent scalar product to the standard $H(\curl;\Omega)$ scalar product.
Lateron we will extend our analysis in Section~\ref{subsec:topnontrivial} to topologically nontrivial domains.
Recall that the variational formulation of \eqref{eq:Steklov-em} is to find $(\lambda,\bu)\in \setC\times H(\curl;\Omega)\setminus\{0\}$ such that
\begin{align}\label{eq:varform}
\spl\bbmu_h^{-1}\curl\bu,\curl\bu'\spr_{\boldL^2(\Omega)}
-\omega^2 \spl\bbespilon_h\bu,\bu'\spr_{\boldL^2(\Omega)}
-\lambda \spl S\bu,S\bu'\spr_{\boldL^2(\partial\Omega)}=0
\end{align}
for all $\bu'\in H(\curl;\Omega)$.
Testing \eqref{eq:varform} with $\nabla w, w\in H^1_*(\Omega)$ we obtain that any solution $(\lambda,\bu)\in \setC\times H(\curl;\Omega)\setminus\{0\}$ satisfies $\bu\in V_h$, and trivially $(\lambda,\bu)$ satisfies \eqref{eq:varform} for all $\bu'\in V_h$.
On the other hand if $(\lambda,\bu)\in \setC\times V_h\setminus\{0\}$ satisfies \eqref{eq:varform} for all $\bu'\in V_h$, then $(\lambda,\bu)\in \setC\times V_h\setminus\{0\}$ satisfies \eqref{eq:varform} also for all $\bu'\in H(\curl;\Omega)$.
Thus we can consider the variational formulation of the eigenvalue problem \eqref{eq:Steklov-em} in the space $V_h$ instead of $H(\curl;\Omega)$.
This is advantageous, because the space $V_h$ admits better embedding properties than $H(\curl;\Omega)$.
However, when we compare now the eigenvalue problems for $h=0$ and $h>0$ we recognize that these problems are posed now on different spaces $V_0\neq V_h, h>0$.
To circumvent the latter we are going to reformulate the eigenvalue problem for $h>0$ in terms of the space $V$.
Note that $\PVh \in L(V,V_h)$.
Since $\curl \PVh \bu= \curl \bu$ for $\bu\in V$ and $\Omega$ is simply connected it follows that $\PVh$ ins injective.
Due to $V_h = \PVh H(\curl;\Omega)$, $H(\curl;\Omega)=V\oplus^\calT \nabla H^1_*(\Omega)$ and $\PVh\nabla H^1_*(\Omega) = \{0\}$ it follows that $\PVh$ surjective.
Thus $\PVh$ is bijection between $V$ and $V_h$.
Consequently we can reformulate the eigenvalue problem for $h>0$ to find $(\lambda,\bu)\in \setC\times V\setminus\{0\}$ such that
\begin{align}\label{eq:varformh}
\spl\bbmu_h^{-1}\curl\bu,\curl\bu'\spr_{\boldL^2(\Omega)}
-\omega^2 \spl\bbespilon_h\PVh\bu,\PVh\bu'\spr_{\boldL^2(\Omega)}
-\lambda \spl S\bu,S\bu'\spr_{\boldL^2(\partial\Omega)}=0
\end{align}
for all $\bu'\in V$, whereat we used that $\curl \PVh \bu= \curl \bu$ and $S\PVh\bu=S\bu$.
We introduce $\Amuh, \Aepsh, \Atr \in L(V)$ through
\begin{align*}
\spl \Amuh\bu,\bu'\spr_V&:=\spl\bbmu_h^{-1}\curl\bu,\curl\bu'\spr_{\boldL^2(\Omega)},\\
\spl \Aepsh\bu,\bu'\spr_V&:=\spl\bbespilon_h\PVh\bu,\PVh\bu'\spr_{\boldL^2(\Omega)}\\
\spl \Atr\bu,\bu'\spr_V&:=\spl S\bu,S\bu'\spr_{\boldL^2(\partial\Omega)}
\end{align*}
for all $\bu,\bu'\in V$.
Note that $A_{\bbmu_0}=I$.
Let $T_h(\lambda):=\Amuh^{-1}(\frac{\omega^2}{\lambda}\Aepsh+\Atr)$ and $\setT_h(\lambda):=I-\lambda T_h(\lambda)$ for $\lambda\in\setC\setminus\{0\}$.
Recall that due to Assumption~\ref{ass:NoDirEV-MW} $\lambda=0$ is not an eigenvalue of \eqref{eq:Steklov-em}.
Thence the eigenvalue problem \eqref{eq:varformh} is equivalent to find $(\lambda,\bu)\in \setC\setminus\{0\}\times V\setminus\{0\}$ such that $\setT_h(\lambda)\bu=0$.
We recall the three-term H\"older's inequality
\begin{align}\label{eq:Hoelder}
\|f_1f_2f_3\|_{L^1(\Omega)} \leq \|f_1\|_{L^p(\Omega)}\|f_2\|_{L^q(\Omega)}\|f_3\|_{L^rp(\Omega)}, \qquad \frac{1}{p}+\frac{1}{q}+\frac{1}{r}=1,
\end{align}
for $p,q,r\in[1,\infty]$.
We also recall the continuous Sobolev embedding for
\begin{align}\label{eq:SobEmbed}
t\in(0,3/2),\, t^*:=6/(3-2t), \quad
\|u\|_{L^{t^*}(\Omega)} \leq C_{L,t} \|u\|_{H^t(\Omega)}
\quad\forall u\in H^t(\Omega),
\end{align}
with the embedding constant $C_{L,t}>0$.
\begin{lemma}\label{lem:norm-est}
\begin{sloppypar}
Let $\Omega$ be a simply connected Lipschitz domain.
Let $(\bbmu_h)_{h\geq0}, (\bbespilon_h)_{h\geq0}$ satisfy Assumption~\ref{ass:mueps}, $\bbmu_0\in\setP\setW^{1,\infty}(\Omega,\calP_1)$, $\bbespilon_0\in\setP\setW^{1,\infty}(\Omega,\calP_2)$, and $s_{\bbmu_0}, s_{\bbespilon_0}>0$ be the maximal Sobolev embedding indices of $H(\curl,\div\bbmu_0,\tr_{\nv\times}^0;\Omega)$ and $H(\curl,\div\bbespilon_0,\tr_{\nv\cdot\bbespilon_0}^0;\Omega)$ (see Proposition~\ref{prop:Ciarlet}) and $\tilde s_{\bbmu_0}:=\min\{1/2,s_{\bbmu_0}\}$.
Let $\tilde s_{\bbmu_0}^*, s_{\bbespilon_0}^*$ be as defined in \eqref{eq:SobEmbed}.
Then for each $p_1\in (\frac{2}{1-2/\tilde s_{\bbmu_0}^*},\infty)$, $p_2\in (\frac{2}{1-2/s_{\bbespilon_0}^*},\infty)$, there exists a constant $C>0$ such that
\begin{align*}
\|T_0(\lambda)-T_h(\lambda)\|_{L(V)} \leq C (1+1/|\lambda|)
\left( \|\bbmu_0-\bbmu_h\|_{\setL^{p_1}(\Omega)}
+ \|\bbespilon_0-\bbespilon_h\|_{\setL^{p_2}(\Omega)} \right)
\end{align*}
for all $\lambda\in\setC\setminus\{0\}$.
The constant $C$ depends only on $\Omega$, $\omega$, $\bbmu_0$, $\bbespilon_0$, $p_1$, $p_2$, and the infinums and supremums in Assumption~\ref{ass:mueps}.
\end{sloppypar}
\end{lemma}
\begin{proof}
We compute
\begin{align}\label{eq:T-Th}
T_0(\lambda)-T_h(\lambda) &= \Amuh^{-1}(\Amuh-I)\left(\frac{\omega^2}{\lambda}\Aepsz+\Atr\right)
+\frac{\omega^2}{\lambda} \Amuh^{-1} (\Aepsz-\Aepsh).
\end{align}
Note that due to Assumption~\ref{ass:mueps} it holds $\sup_{h\geq0}\|\Amuh^{-1}\|_{L(V)}<\infty$.
First we estimate $\|\Aepsz-\Aepsh\|_{L(V)}$.
We compute for $\bu,\bu'\in V$ that
\begin{align}\label{eq:EstEps1}
\begin{aligned}
\spl (\Aepsz-\Aepsh)\bu,\bu' \spr_{L(V)}&=
\spl \bbespilon_0 \bu, \bu' \spr_{\boldL^2(\Omega)}
-\spl \bbespilon_h \PVh \bu, \PVh \bu' \spr_{\boldL^2(\Omega)}\\
&=\spl \bbespilon_0 \bu, \bu' \spr_{\boldL^2(\Omega)}
-\spl \bbespilon_h \PVh \bu, \bu' \spr_{\boldL^2(\Omega)}\\
&=\spl (\bbespilon-\bbespilon_h) \bu, \bu' \spr_{\boldL^2(\Omega)}
+\spl \bbespilon_h (1-\PVh)\bu, \bu' \spr_{\boldL^2(\Omega)}.
\end{aligned}
\end{align}
Let $t_2^*:=2/(1-2/p_2)$, i.e.\ $1/p_2+1/t_2^*+1/2=1$.
Let $t_2:=3/2-3/t_2^*$, i.e.\ $t_2^*=6/(3-2t_2)$.
Let $C_{H,t}>$ be the $V\hookrightarrow\boldH^t(\Omega)$ embedding constant.
Recall that $(1-\PVh)\boldu=\nabla w$ with $w\in H^1_*(\Omega)$ being the solution to~\eqref{eq:eqforw}, and that $\epsilon_-$ is defined in Assumption~\ref{ass:mueps}.
It follows 
\begin{align*}
\|\nabla w\|_{\boldL^2(\Omega)} \leq \epsilon_-^{-1} \|\spl \bbespilon_h \bu, \nabla \cdot\spr \|_{H^{-1}(\Omega)}
&=\epsilon_-^{-1} \|\spl (\bbespilon_0-\bbespilon_h) \bu, \nabla \cdot\spr \|_{H^{-1}}\\
&\leq \epsilon_-^{-1} \|\bbespilon_0-\bbespilon_h\|_{\setL^{p_2}(\Omega)} \|\bu\|_{\boldL^{t_2^*}(\Omega)}\\
&\leq C_{L,t_2} \epsilon_-^{-1} \|\bbespilon_0-\bbespilon_h\|_{\setL^{p_2}(\Omega)} \|\bu\|_{\boldH^{t_2}(\Omega)}\\
&\leq C_{H,t_2} C_{L,t_2} \epsilon_-^{-1} \|\bbespilon_0-\bbespilon_h\|_{\setL^{p_2}(\Omega)} \|\bu\|_{V}.
\end{align*}
Thence
\begin{align}\label{eq:EstEps2}
|\spl \bbespilon_h (1-&\PVh)\bu, \bu' \spr_{\boldL^2(\Omega)}|
\leq (\sup_{h\geq0}\|\bbespilon_h\|_{\setL^\infty(\Omega)}) C_{H,0}
C_{H,t_2} C_{L,t_2} \epsilon_-^{-1} \|\bbespilon_0-\bbespilon_h\|_{\setL^{p_2}(\Omega)} \|\bu\|_{V} \|\bu'\|_{V}.
\end{align}
Further we estimate
\begin{align}\label{eq:EstEps3}
\begin{aligned}
|\spl (\bbespilon_0-\bbespilon_h) \bu, \bu' \spr_{\boldL^2(\Omega)}|
&\leq \|\bbespilon_0-\bbespilon_h\|_{\setL^{p_2}(\Omega)} \|\bu\|_{\boldL^{t_2^*}(\Omega)} \|\bu'\|_{\boldL^{2}(\Omega)} \\
&\leq C_{L,t_2} C_{H,0} \|\bbespilon_0-\bbespilon_h\|_{\setL^{p_2}(\Omega)} \|\bu\|_{\boldH^{t_2}(\Omega)} \|\bu'\|_{V}\\
&\leq C_{H,t_2} C_{L,t_2} C_{H,0} \|\bbespilon_0-\bbespilon_h\|_{\setL^{p_2}(\Omega)} \|\bu\|_{V} \|\bu'\|_{V}.
\end{aligned}
\end{align}
Together \eqref{eq:EstEps1}, \eqref{eq:EstEps2},\eqref{eq:EstEps3} yield that
\begin{align}\label{eq:EstEpsTot}
\|\Aepsz-\Aepsh\|_{L(V)}\leq C \|\bbespilon_0-\bbespilon_h\|_{\setL^{p_2}(\Omega)}
\end{align}
with a constant $C>0$.\\
Next we estimate $\|(\Amuh-I)\Aepsz\|_{L(V)}$.
For $\bu\in V$ let $\bv_1:=\Aepsz\bu$.
Thence $\bv_1$ satisfies
\begin{align*}
\spl\bbmu_0^{-1}\curl\bv_1,\curl\bu'\spr_{\boldL^2(\Omega)}
=\spl\bv_1,\bu'\spr_V
=\spl\Aepsz\bu,\bu'\spr_V
=\spl \bbespilon_0\bu, \bu'\spr_{\boldL^2(\Omega)}
\end{align*}
for all $\bu'\in V$, and even also for all $\bu'\in H(\curl;\Omega)$.
Hence $\curl\bbmu_0^{-1}\curl\bv_1=\bbespilon_0\bu\in\boldL^2(\Omega)$ and $\nv\times\bbmu_0^{-1}\bv_1=0$ on $\partial\Omega$, i.e.\ $\bbmu_0^{-1}\curl\bv_1 \in H(\curl,\div\bbmu_0,\tr_{\nv\times}^0;\Omega)$.
Let $t_1^*:=2/(1-2/p_1)$, i.e.\ $1/p_1+1/t_1^*+1/2=1$.
Let $t_1:=3/2-3/t_1^*$, i.e.\ $t_1^*=6/(3-2t_1)$.
Recall that $\mu_-$ is defined in Assumption~\ref{ass:mueps}, and let $\tilde C_{H,t}$ be the $H(\curl,\div\bbmu_0,\tr_{\nv\times}^0;\Omega)\hookrightarrow\boldH^t(\Omega)$ embedding constant.
We estimate for $\bu,\bu'\in V$ that
\begin{align}\label{eq:EstAmu1}
\begin{aligned}
|\spl(\Amuh-I)\Aepsz\bu,\bu'\spr_{L(V)}|
&=|\spl (\bbmu_h^{-1}-\bbmu_0^{-1}) \curl\bv_1,\curl\bu'\spr_{\boldL^2(\Omega)}|\\
&=|\spl \bbmu_h^{-1}(\bbmu_0-\bbmu_h)\bbmu_0^{-1} \curl\bv_1,\curl\bu'\spr_{\boldL^2(\Omega)}|\\
&\leq (\sup_{h\geq0}\|\bbmu_h^{-1}\|_{\setL^\infty(\Omega)}) \|\bbmu_0-\bbmu_h\|_{\setL^{p_1}(\Omega)} \|\bbmu_0^{-1} \curl\bv_1\|_{\boldL^{t_1^*}(\Omega)} \|\curl\bu'\|_{\boldL^2(\Omega)}\\
&\leq C_{L,t_1} \mu_-^{-3/2} \|\bbmu_0-\bbmu_h\|_{\setL^{p_1}(\Omega)} \|\bbmu_0^{-1}\curl\bv_1\|_{\boldH^t(\Omega)} \|\bu'\|_{V}\\
&\leq \tilde C_{H,t_1} C_{L,t_1} \mu_-^{-3/2} \|\bbmu_0-\bbmu_h\|_{\setL^{p_1}(\Omega)} \|\bbmu_0^{-1}\curl\bv_1\|_{H(\curl;\Omega)} \|\bu'\|_{V}\\
&\leq C \|\bbmu_0-\bbmu_h\|_{\setL^{p_1}(\Omega)} \|\bu\|_{V} \|\bu'\|_{V},
\end{aligned}
\end{align}
with $C=(C_{H,0}\mu_-^{-1/2}+1) C_{H,0} \|\bbespilon_0\|_{\setL^\infty(\Omega)} \sqrt{2}
\tilde C_{H,t_1} C_{L,t_1} \mu_-^{-3/2}$,
whereat we used that
\begin{align*}
\|\bbmu_0^{-1}\curl\bv_1\|_{H(\curl;\Omega)}
&\leq \sqrt{2} (\|\bbmu_0^{-1}\curl\bv_1\|_{\boldL^2(\Omega)}+\|\bbespilon_0\bu\|_{\boldL^{2}(\Omega)})\\
&\leq \mu_-^{-1/2}\sqrt{2} \|\bv_1\|_V
+C_{H,0} \|\bbespilon_0\|_{\setL^\infty(\Omega)} \sqrt{2} \|\bu\|_V\\
&\leq C_{H,0}^2 \|\bbespilon_0\|_{\setL^\infty(\Omega)} \mu_-^{-1/2}\sqrt{2} \|\bu\|_V
+C_{H,0} \|\bbespilon_0\|_{\setL^\infty(\Omega)} \sqrt{2} \|\bu\|_V.
\end{align*}
Hence
\begin{align}\label{eq:EstAmu2}
\|(\Amuh-I)\Aepsz\|_{L(V)}\leq C\|\bbmu_0-\bbmu_h\|_{\setL^{p_1}(\Omega)}.
\end{align}
Next we estimate $\|(\Amuh-I)\Atr\|_{L(V)}$.
For $\bu\in V$ let $\bv_2:=\Atr\bu$.
In a similar way as for $\bv_1$ we obtain that $\curl\bbmu_0^{-1}\curl\bv_2=0$ and $\nv\times\bbmu_0^{-1}\curl\bv_2=S\bu\in\boldL^2(\partial\Omega)$.
Since $\Omega$ is simply connected and $\curl\bbmu_0^{-1}\curl\bv_2=0$ it follows that $\bbmu_0^{-1}\curl\bv_2=\nabla w$ for a $w\in H^1(\Omega)$.
Thence $\nv\times\bbmu_0^{-1}\curl\bv_2\in\boldL^2(\partial\Omega)$ implies that $\tr w\in H^1(\partial\Omega)$.
Note that this point of the proof requires us to work with $\tilde s_\bbmu$ instead of $s_\bbmu$.
It also holds that $\div\bbmu\nabla w=\div\curl\bv_2=0$ in $\Omega$.
Let $R$ be a bounded extension operator $H^1(\partial\Omega)\to H^{1+t_1}(\Omega)$.
Then $w_1:=w-R\tr w\in H^1_0(\Omega)$ solves
\begin{align*}
\spl\bbmu_0\nabla w_1,\nabla w' \spr_{\boldL^2(\Omega)}
=-\spl\bbmu_0 \nabla R\tr w,\nabla w' \spr_{\boldL^2(\Omega)}
\end{align*}
for all $w'\in H^1_0(\Omega)$.
Note that the regularity index $s_{\bbmu_0}$ equals the regularity index for the former scalar Dirichlet problem (see the proof of Proposition~\ref{prop:Ciarlet}).
The right hand-side of the former equation is in $H^{-1+t_1}(\Omega)$.
Thus it follows with \cite[Theorem~6.8]{Ciarlet:20} that
\begin{align*}
\|\bbmu_0^{-1}\curl\bv_2\|_{\boldH^{t_1}(\Omega)}\leq \|\nabla w_1\|_{\boldH^{t_1}(\Omega)}+\|\nabla R\tr w\|_{\boldH^{t_1}(\Omega)} \leq C \|\bu\|_V
\end{align*}
with a constant $C>0$.
Hence estimates similar to \eqref{eq:EstAmu1} yield that
\begin{align}\label{eq:EstAmu3}
\|(\Amuh-I)\Atr\|_{L(V)}\leq C\|\bbmu-\bbmu_h\|_{\setL^{p_1}(\Omega)}
\end{align}
with a constant $C>0$.
The claim follows now from \eqref{eq:T-Th}, \eqref{eq:EstEpsTot}, \eqref{eq:EstAmu2} and \eqref{eq:EstAmu3}.
\end{proof}
\begin{theorem}\label{thm:ev-pert}
\begin{sloppypar}
Let $\Omega$ be a simply connected Lipschitz domain.
Let $(\bbmu_h)_{h\geq0}, (\bbespilon_h)_{h\geq0}$ satisfy Assumption~\ref{ass:mueps}, $\bbmu_0\in\setP\setW^{1,\infty}(\Omega,\calP_1)$, $\bbespilon_0\in\setP\setW^{1,\infty}(\Omega,\calP_2)$, and $s_{\bbmu_0}, s_{\bbespilon_0}>0$ be the maximal Sobolev embedding indices of $H(\curl,\div\bbmu_0,\tr_{\nv\times}^0;\Omega)$ and $H(\curl,\div\bbespilon_0,\tr_{\nv\cdot\bbespilon_0}^0;\Omega)$ (see Proposition~\ref{prop:Ciarlet}) and $\tilde s_{\bbmu_0}:=\min\{1/2,s_{\bbmu_0}\}$.
Let $\tilde s_{\bbmu_0}^*, s_{\bbespilon_0}^*$ be as defined in \eqref{eq:SobEmbed} and $p_1\in (\frac{2}{1-2/\tilde s_{\bbmu_0}^*},\infty)$, $p_2\in (\frac{2}{1-2/s_{\bbespilon_0}^*},\infty)$.
Let $\lim_{h\to0}\|\bbmu_0-\bbmu_h\|_{\setL^{p_1}(\Omega)}=\lim_{h\to0}\|\bbespilon_0-\bbespilon_h\|_{\setL^{p_2}(\Omega)}=0$.
Then
\begin{enumerate}
 \item if $(\lambda_h)_{h>0}$ is a sequence of eigenvalues of \eqref{eq:Steklov-em} with $(\bbmu,\bbespilon)=(\bbmu_h,\bbespilon_h)$, which converges to $\lambda_0$, then $\lambda_0$ is an eigenvalue of \eqref{eq:Steklov-em} with $(\bbmu,\bbespilon)=(\bbmu_0,\bbespilon_0)$;
 \item for every eigenvalue $\lambda_0$ of \eqref{eq:Steklov-em} with $(\bbmu,\bbespilon)=(\bbmu_0,\bbespilon_0)$ exists a sequence $(\lambda_h)_{h>0}$ with $\lambda_h$ being an eigenvalue of \eqref{eq:Steklov-em} with $(\bbmu,\bbespilon)=(\bbmu_h,\bbespilon_h)$ for almost all $h>0$ such that $(\lambda_h)_{h>0}$ converges to $\lambda_0$;
 \item let $\Lambda\subset\setC$ be compact and such that $\partial\Lambda$ is contains no eigenvalues of \eqref{eq:Steklov-em} with $(\bbmu,\bbespilon)=(\bbmu_0,\bbespilon_0)$. Then there exists an index $h_0>0$ such that for all $h<h_0$ the sum of the algebraic multiplicities of all eigenvalues of \eqref{eq:Steklov-em} in $\Lambda$ are equal for $(\bbmu,\bbespilon)=(\bbmu_0,\bbespilon_0)$ and $(\bbmu,\bbespilon)=(\bbmu_h,\bbespilon_h)$;
 \item let $\lambda_0\in\setC$ be an eigenvalue of \eqref{eq:Steklov-em} with $(\bbmu,\bbespilon)=(\bbmu_0,\bbespilon_0)$ and $\kappa$ be the maximal length of its Jordan chains. Let $\Lambda\subset\setC$ be compact, such that the $\lambda_0\in\Lambda$ is the only eigenvalue of \eqref{eq:Steklov-em} with $(\bbmu,\bbespilon)=(\bbmu_0,\bbespilon_0)$ contained in the closure of $\Lambda$.
 Let $\lambda_h^\mathrm{mean}$ be the weighted mean (with respect to the algebraic multiplicity) of all eigenvalues of \eqref{eq:Steklov-em} with $(\bbmu,\bbespilon)=(\bbmu_h,\bbespilon_h)$ in $\Lambda$.
 Let $\lambda_h\in\Lambda$ be an eigenvalue of \eqref{eq:Steklov-em} with $(\bbmu,\bbespilon)=(\bbmu_h,\bbespilon_h)$. Then it holds that
 \begin{align*}
 |\lambda_0-\lambda_h|=O\left(
 \|\bbmu_0-\bbmu_h\|_{\setL^{p_1}(\Omega)}^\frac{1}{\kappa}
 +\|\bbespilon_0-\bbespilon_h\|_{\setL^{p_2}(\Omega)}^\frac{1}{\kappa}
 \right)
 \end{align*}
 and
 \begin{align*}
 |\lambda_0-\lambda_h^\mathrm{mean}|=O\left(
 \|\bbmu_0-\bbmu_h\|_{\setL^{p_1}(\Omega)}
 +\|\bbespilon_0-\bbespilon_h\|_{\setL^{p_2}(\Omega)}
 \right).
 \end{align*}
\end{enumerate}
\end{sloppypar}
\end{theorem}
\begin{proof}
We apply \cite{Karma:96a,Karma:96b} to $\setT_0(\cdot)$, $\setT_h(\cdot)$ and use Lemma~\ref{lem:norm-est}.
The first two claims follow from \cite[Theorem~2]{Karma:96a}.
The third claim follows from \cite[Theorem~3]{Karma:96a}.
The fourth claim follows from \cite[Theorem~2]{Karma:96b}.
\end{proof}
Note that in the following Theorem~\ref{thm:improved-rate} the allowed range of $p_1$ is $p_1\in(\frac{1}{1-2/\tilde s_{\bbmu_0}^*},\infty)$, whereas in Theorem~\ref{thm:ev-pert} it is required that $p_1\in(\frac{2}{1-2/\tilde s_{\bbmu_0}^*},\infty)$.
Further note that \eqref{eq:est-improved} does not depend on $\kappa$.
\begin{theorem}\label{thm:improved-rate}
Let the assumptions of Theorem~\ref{thm:ev-pert} be satisfied.
Let $\lambda_0$ be an eigenvalue of \eqref{eq:Steklov-em} with $(\bbmu,\bbespilon)=(\bbmu_0,\bbespilon_0)$ and $(\bu_n), n=1,\dots,N$ be a normalized basis of the corresponding generalized eigenspace.
Assume that
\begin{align}\label{eq:neqz}
1+\frac{\lambda_0^2}{N}\sum_{n=1}^N \spl T_0'(\lambda_0)\bu_n,\ol{\bu_n}\spr_{V}\neq0,
\end{align}
whereat $T_0'(\lambda)$ is the derivative of $T_0(\lambda)$ with respect to $\lambda$ and $\ol{\bu_n}$ denotes the complex conjugation of $\bu_n$.
Then the last statement of Theorem~\ref{thm:ev-pert} improves with $p_1\in (\frac{1}{1-2/\tilde s_{\bbmu_0}^*},\infty)$ to
\begin{align}\label{eq:est-improved}
|\lambda_0-\lambda_h|=O\left(
\|\bbmu_0-\bbmu_h\|_{\setL^{p_1}(\Omega)}
+\|\bbespilon_0-\bbespilon_h\|_{\setL^{p_2}(\Omega)}
\right).
\end{align}
\end{theorem}
\begin{proof}
First note that $(\ol{\bu_n}), n=1,\dots,N$ is a normalized basis of the generalized eigenspace of the adjoint problem.
Also note that $2p_1\in (\frac{2}{1-2/\tilde s_{\bbmu_0}^*},\infty)$.
We apply \cite[Theorem~4.1]{Moskow:15} to $T_0(\cdot), T_h(\cdot)$ and use Lemma~\ref{lem:norm-est}.
Let $C\neq0$ be the left hand-side of \eqref{eq:neqz}.
Thence
\begin{align*}
\frac{CN}{\lambda_0^2}(\lambda_0-\lambda_h)&= \sum_{n=1}^N
\spl (T_0(\lambda_0)-T_h(\lambda_0))\bu_n,\ol{\bu_n}\spr_V\\
&+O\left(\|\bbmu_0-\bbmu_h\|^{2}_{\setL^{2p_1}(\Omega)}+\|\bbespilon_0-\bbespilon_h\|^{2}_{\setL^{p_2}(\Omega)}\right).
\end{align*}
Recall~\eqref{eq:T-Th}:
\begin{align*}
T_0(\lambda_0)-T_h(\lambda_0) &= \Amuh^{-1}(\Amuh-I)\left(\frac{\omega^2}{\lambda_0}\Aepsz+\Atr\right)
+\frac{\omega^2}{\lambda_0} \Amuh^{-1} (\Aepsz-\Aepsh).
\end{align*}
We estimate as in the proof of Lemma~\ref{lem:norm-est} that
\begin{align*}
\spl \Amuh^{-1} (\Aepsz-\Aepsh)\bu_n,\ol{\bu_n}\spr_{V}=
O\left(\|\bbespilon_0-\bbespilon_h\|_{\setL^{p_2}(\Omega)}\right).
\end{align*}
For the other term we compute
\begin{align*}
\Amuh^{-1} (\Amuh-I)\left(\frac{\omega^2}{\lambda_0}\Aepsz+\Atr\right)
&=(\Amuh-I)\left(\frac{\omega^2}{\lambda_0}\Aepsz+\Atr\right)\\
&-(\Amuh-I)\Amuh^{-1}(\Amuh-I)\left(\frac{\omega^2}{\lambda_0}\Aepsz+\Atr\right)
\end{align*}
and estimate with $\bv_n:=(\frac{\omega^2}{\lambda_0}\Aepsz+\Atr)\bu_n$ that
\begin{align*}
\Big|\big\spl \Amuh^{-1} (\Amuh-I)
\Big(\frac{\omega^2}{\lambda_0}\Aepsz+\Atr\Big)\bu_n,&\ol{\bu_n}
\big\spr_V\Big|\\
&\leq \|\bbmu\|_{\setL^\infty(\Omega)} \mu_-^{-1} |\spl (\bbmu_h-\bbmu_0) \bbmu_0^{-1}\curl\bv_n,\bbmu_0^{-1}\curl\ol{\bu_n}\spr_{\boldL^2(\Omega)}|\\
&+(\sup_{h\geq0}\|\bbmu_h^{-1}\|_{\boldL^\infty(\Omega)}) \mu_-^{-1} \|(\Amuh-I)\bv_n\|_V \|(\Amuh-I)\ol{\bu_n}\|_V.
\end{align*}
Hence by means of \eqref{eq:Hoelder} we estimate similarly as in the proof of Lemma~\ref{lem:norm-est} further that
\begin{align*}
|\spl (\bbmu_h-\bbmu_0) \bbmu_0^{-1}\curl\bv_n,\bbmu_0^{-1}\curl\ol{\bu_n}\spr_{\boldL^2(\Omega)}|
=O\left( \|\bbmu_0-\bbmu_h\|_{\setL^{p_1}(\Omega)} \right)
\end{align*}
and
\begin{align*}
\|(\Amuh-I)\bv_n\|_V \|(\Amuh-I)\ol{\bu_n}\|_V
=O\left( \|\bbmu_0-\bbmu_h\|^2_{\setL^{2p_1}(\Omega)} \right),
\end{align*}
whereat we exploited that $\bbmu_0^{-1}\curl\ol{\bu_n}$ has the same regularity properties as $\bbmu_0^{-1}\curl\bv_n$.
We recall H\"older's interpolation inequality
\begin{align}\label{eq:HoelderInt}
\|f\|_{L^r} \leq \|f\|_{L^p}^\alpha \|f\|_{L^q}^{1-\alpha}
\end{align}
for $q\in[1,\infty]$, $p\in[1,q]$, $r\in[p,q]$ and $\alpha\in[0,1]$ with $\frac{1}{r}=\frac{\alpha}{p}+\frac{1-\alpha}{q}$.
At last we estimate as in \cite{Cogar:22} by means of \eqref{eq:HoelderInt} with $q=\infty$, $p=p_1$, $r=2p_1$, $\alpha=1/2$ that
\begin{align*}
\|\bbmu_0-\bbmu_h\|^2_{\setL^{2p_1}(\Omega)}
\leq \|\bbmu_0-\bbmu_h\|_{\setL^{p_1}(\Omega)}.
\end{align*}
Hence altogether we obtain that
\begin{align*}
|\lambda_0-\lambda_h|=O\left(
\|\bbmu_0-\bbmu_h\|_{\setL^{p_1}(\Omega)}
+\|\bbespilon_0-\bbespilon_h\|_{\setL^{p_2}(\Omega)}
\right).
\end{align*}
\end{proof}
\begin{remark}
Note that for a semisimple eigenvalue (i.e.\ $\kappa=1$) the left hand-side of \eqref{eq:neqz} equals $\frac{1}{N}\sum_{n=1}^N \spl S\bu_n,S\ol{\bu_n}\spr_{\boldL^2(\partial\Omega)}$.
If we compare Theorem~\ref{thm:improved-rate} with $N=1$ to \cite[Theorem~7]{Cogar:22}, we observe that here we require \eqref{eq:neqz} whereas no such assumptions occurs in \cite[Theorem~7]{Cogar:22}.
The explanation is that there is a slight mistake in \cite{Cogar:22}.
Indeed \cite[Theorem~6]{Cogar:22} is an incorrect version of \cite[Theorem~4.1]{Moskow:15}, which assumes that the eigenspaces for $T_0(\cdot)$ and its adjoint problem are equal.
However, in general this is only true for selfadjoint operators.
Thus in the right arguments of scalar products in \cite[Theorems~6,7]{Cogar:22} $u_0$ needs to be replaced by $u_0^*$.
For Steklov eigenvalue problems it holds $u_0^*=\ol{u_0}$ (with the overline indicating the complex conjugation).
Thus one cannot deduce from $S\bu\neq0$ that $\spl S\bu,S\ol{\bu}\spr_{\boldL^2(\partial\Omega)}\neq0$.
\end{remark}
\begin{remark}
Note that for Theorem~\ref{thm:improved-rate} we apply \cite[Theorem~4.1]{Moskow:15} and can obtain an explicit formula for the leading order perturbation term with respect to perturbations in $\bbepsilon$.
However, for perturbations in $\bbmu$ the explicit perturbation term turns out to be of the same order as the remainder.
\end{remark}

\subsection{Topologically nontrivial domains}\label{subsec:topnontrivial}

Until now we have assumed in this section that the domain $\Omega$ is simply connected.
Let us discuss now how to extend the former analysis and in particular Theorems~\ref{thm:ev-pert} and~\ref{thm:improved-rate} to domains $\Omega$ which satisfy Assumption~\ref{ass:top} with $\mathbf{(Top)}_{I>0}$.
To construct a suitable topological decomposition of $H(\curl;\Omega)$ in this case we follow \cite[Chapter~3]{AssousCiarletLabrunie:18} with some minor adaptations.
Let
\begin{align*}
P_*(\dot\Omega)&:=\{\dot w\in H^1_*(\dot\Omega)\colon [\dot w]_{\Sigma_i}=\setC\quad\text{for all}\quad i=1,\dots,I\},
\end{align*}
whereat $[\dot w]_{\Sigma_i}$ denotes the jump of $w$ on $\Sigma_i$.
Functions $\dot w\in P_*(\dot\Omega)$ in general do not admit a weak gradient $\nabla \dot w$ in $\Omega$.
However, the weak gradient $\nabla \dot w$ is well defined in $L^2(\dot\Omega)$, which we can identify with $L^2(\Omega)$.
We use the notation $\widetilde{\nabla}\dot w\in L^2(\Omega)$ for the former function to emphasize the difference.
Note that $\curl\widetilde{\nabla}\dot w=0$ and $S\widetilde{\nabla}\dot w=0$ for $\dot w\in P_*(\dot\Omega)$.
Recall \cite[Theorem~3.3.2]{AssousCiarletLabrunie:18} which tells that for $\bu\in H(\curl;\Omega)$, $\curl\bu=0$, if and only if $\bu=\widetilde{\nabla}\dot w$ with $\dot w\in P_*(\dot\Omega)$.
Let
\begin{align*}
\tilde V_h&:=\{\bu\in H(\curl;\Omega)\colon \spl \bbespilon_h\bu,\widetilde{\nabla}\dot w\spr_{\boldL^2(\Omega)}=0 \quad\text{for all}\quad \dot w\in P_*(\dot\Omega)\}.
\end{align*}
The projection onto $\tilde V_h$ is given by $\PVht\bu:=\bu-\widetilde{\nabla}\dot w$ with $\dot w\in P_*(\dot\Omega)$ being the solution to
\begin{align}\label{eq:eqforwt}
\spl \bbespilon_h \widetilde{\nabla}\dot w, \widetilde{\nabla}\dot w' \spr_{\boldL^2(\Omega)} = \spl \bbespilon_h \bu, \widetilde{\nabla}\dot w' \spr_{\boldL^2(\Omega)} \quad\text{for all } \dot w'\in P_*(\dot\Omega).
\end{align}
\begin{sloppypar}
\noindent
Hence we have a topological decomposition $H(\curl;\Omega)=\tilde V_h\oplus^\calT \widetilde{\nabla}P_*(\dot\Omega)$ and $\spl\bu,\bu'\spr_{\tilde V_h}$ $:=$ $\spl \bbmu_h^{-1}\curl\bu,\curl\bu'\spr_{\boldL^2(\Omega)}$ is on $\tilde V_h$ an equivalent scalar product to the standard $H(\curl;\Omega)$ scalar product.
Thus we can repeat the preceding analysis to reformulate eigenvalue problem \eqref{eq:Steklov-em} in terms of $\tilde\setT_h(\lambda)=\setI-\lambda \tilde T_h(\lambda)$ on the space $\tilde V:=\tilde V_0$, whereat in the definitions of $\tilde T_h(\cdot)$ quantities $V, \PVh$ merely need to be replaced by $\tilde V, \tilde \PVh$, etc..
\end{sloppypar}
\begin{theorem}\label{thm:topnontrivial}
Let $\Omega$ satisfy Assumption~\ref{ass:top} with $\mathbf{(Top)}_{I>0}$.
Then Theorems~\ref{thm:ev-pert} and \ref{thm:improved-rate} still hold true, whereat in the latter Theorem~\ref{thm:improved-rate} quantities $V, T_0(\cdot)$ need to be replaced by $\tilde V, \tilde T_0(\cdot)$.
\end{theorem}
\begin{proof}
First we show that Lemma~\ref{lem:norm-est} can be generalized.
Note that Proposition~\ref{prop:Ciarlet} holds already for $\mathbf{(Top)}_{I>0}$ and hence Estimates~\eqref{eq:EstEpsTot}, \eqref{eq:EstAmu2} do not require any adaptation.
To obtain \eqref{eq:EstAmu3} it suffices to bound $\|\bbmu_0^{-1}\curl\bv_2\|_{\boldH^{t_1}(\Omega)}$, $\bv_2:=\Atr\bu$ in terms of $\|\bu\|_{\tilde V}$.
It holds $\curl \bbmu_0^{-1}\curl\bv_2=0$ and $\nv\times\bbmu_0^{-1}\curl\bv_2=S\bu\in \boldL^2(\partial\Omega)$.
Due to \cite[Theorem~3.3.2]{AssousCiarletLabrunie:18} we obtain $\bbmu_0^{-1}\curl\bv_2=\widetilde{\nabla}\dot w$ for a $\dot w\in P_*(\dot\Omega)$.
Note that functions in $P_*(\dot\Omega)$ have a well-defined trace in $H^s(\partial\Omega)$, $s<1/2$.
From $\nv\times\bbmu_0^{-1}\curl\bv_2=S\bu\in \boldL^2(\partial\Omega)$ we deduce that $\tr\dot w\in H^1(\partial\Omega)$.
Let $w\in H^1(\Omega)$ be the solution to $\div\mu\nabla w=0$ in $\Omega$, $w=\dot w$ on $\partial\Omega$.
As in the proof of Lemma~\ref{lem:norm-est} we can bound $\|\nabla w\|_{\boldH^{t_1}(\Omega)}$ in terms of $\|\bu\|_{\tilde V}$.
In addition it holds $\widetilde{\nabla}\dot w-\nabla w\in H(\curl,\div\bbmu_0,\tr_{\nv\times}^0;\Omega)$ and hence by means of Proposition~\ref{prop:Ciarlet} we can estimate
\begin{align*}
\|\bbmu_0^{-1}\curl\bv_2\|_{\boldH^{t_1}(\Omega)}
&\leq \|\widetilde{\nabla}\dot w-\nabla w\|_{\boldH^{t_1}(\Omega)}+\|\nabla w\|_{\boldH^{t_1}(\Omega)} \\
&\leq \tilde C_{H,t_1} \|\widetilde{\nabla}\dot w-\nabla w\|_{\boldL^2(\Omega)}
+\|\nabla w\|_{\boldH^{t_1}(\Omega)} \\
&\leq \tilde C_{H,t_1} \|\widetilde{\nabla}\dot w\|_{\boldL^2(\Omega)}
+(\tilde C_{H,t_1}+1)\|\nabla w\|_{\boldH^{t_1}(\Omega)} \\
&\leq C_1( \|\bbmu_0^{-1}\curl\bv_2\|_{\boldL^2(\Omega)}
+\|S\bu\|_{\boldL^{2}(\partial\Omega)} )\\
&\leq C_2 \|\bu\|_{\tilde V}
\end{align*}
with a constants $C_1,C_2>0$.
Altogether we established the generalization of Lemma~\ref{lem:norm-est} from which the generalization of Theorem~\ref{thm:ev-pert} follows directly.
For the generalization of Theorem~\ref{thm:improved-rate} we make the same modifications as for the generalization of Lemma~\ref{lem:norm-est}.
\end{proof}

%

\bibliographystyle{abbrv}
\bibliography{short_biblio}
\end{document}

%% file: commands.tex
\usepackage[utf8]{inputenc}
\usepackage[ngerman, english]{babel}
\usepackage{amsmath,amsthm,amssymb,amsfonts}
\usepackage{mathrsfs}
\usepackage{ifthen}
\usepackage{enumerate}
\usepackage{comment}
\usepackage{subcaption}
\usepackage{hyphsubst}
\usepackage{xcolor}
\usepackage{calc}
\usepackage{graphicx}
\usepackage[toc,page]{appendix}
\usepackage[bbgreekl]{mathbbol}
\DeclareSymbolFontAlphabet{\mathbbm}{bbold}
\DeclareSymbolFontAlphabet{\mathbb}{AMSb}
\newcommand{\bbepsilon}{\bbespilon}

\renewcommand{\Re}{\operatorname{Re}}

\newcommand{\bu}{\boldu}
\newcommand{\bv}{\boldv}

\makeatletter
\def\@hspace#1{\begingroup\setlength\dimen@{#1}\hskip\dimen@\endgroup}
\makeatother

\newtheorem{theorem}{Theorem}[section]
\newtheorem{lemma}[theorem]{Lemma}

\newtheorem{remark}[theorem]{Remark}

\newtheorem{assumption}[theorem]{Assumption}
\newtheorem{proposition}[theorem]{Proposition}

\newcommand{\PVh}{P_{V_h}}
\newcommand{\PVht}{P_{\tilde V_h}}
\newcommand{\SC}{K}
\newcommand{\indSC}{s}

\newcommand{\Amuh}{A_{\bbmu_h}}
\newcommand{\Aepsh}{A_{\bbespilon_h}}
\newcommand{\Aepsz}{A_{\bbespilon_0}}
\newcommand{\Atr}{A_\mathrm{tr}}

\newcommand{\nv}{\boldnu}

\newcommand{\ol}[1]{\overline{#1}}

\newcommand{\spl}{\langle}
\newcommand{\spr}{\rangle}
\newcommand{\bpm}{\begin{pmatrix}}
\newcommand{\epm}{\end{pmatrix}}

\DeclareMathOperator{\curl}{curl}

\renewcommand{\div}{\operatorname{div}}

\DeclareMathOperator{\tr}{tr}
\DeclareMathOperator{\ran}{ran}

\newcommand{\boldnu}{\boldsymbol{\nu}}

\newcommand{\setC}{\mathbb{C}}

\newcommand{\setI}{\mathbb{I}}

\newcommand{\setL}{\mathbb{L}}

\newcommand{\setN}{\mathbb{N}}

\newcommand{\setP}{\mathbb{P}}

\newcommand{\setR}{\mathbb{R}}

\newcommand{\setT}{\mathbb{T}}

\newcommand{\setW}{\mathbb{W}}
\newcommand{\setX}{\mathbb{X}}

\newcommand{\boldu}{\mathbf{u}}
\newcommand{\boldv}{\mathbf{v}}

\newcommand{\boldC}{\mathbf{C}}

\newcommand{\boldH}{\mathbf{H}}

\newcommand{\boldL}{\mathbf{L}}

\newcommand{\boldX}{\mathbf{X}}

\newcommand{\calH}{\mathcal{H}}

\newcommand{\calP}{\mathcal{P}}

\newcommand{\calT}{\mathcal{T}}

\definecolor{brickred}{rgb}{0.8, 0.25, 0.33}
\definecolor{bostonuniversityred}{rgb}{0.8, 0.0, 0.0}
\definecolor{cornellred}{rgb}{0.7, 0.11, 0.11}
\definecolor{corn}{rgb}{0.98, 0.93, 0.36}
\definecolor{schoolbusyellow}{rgb}{1.0, 0.85, 0.0}
\definecolor{TUblue}{rgb}{0,102,153}
\colorlet{TUbluelight}{TUblue!30!white}